\documentclass[12pt]{article}
\usepackage{amsfonts}
\usepackage{amssymb}
\usepackage{color}
\usepackage{amsmath}
=0pt
\def\sqr#1#2{{\vcenter{\vbox{\hrule height.#2pt
              \hbox{\vrule width.#2pt height#1pt \kern#1pt \vrule
width.#2pt}
              \hrule height.#2pt}}}}
\def\signed #1{{\unskip\nobreak\hfil\penalty50
              \hskip2em\hbox{}\nobreak\hfil#1
              \parfillskip=0pt \finalhyphendemerits=0 \par}}
\def\endpf{\signed {$\sqr69$}}
\def\dbC{{\mathbb{C}}}

\def\dbN{{\mathbb{N}}}

\def\dbP{{\mathbb{P}}}

\def\dbR{{\mathbb{R}}}

\def\a{\alpha}

\def\g{\gamma}
\def\d{\delta}
\def\e{\varepsilon}

\def\k{\kappa}
\def\l{\lambda}

\def\f{\varphi}

\def\o{\omega}

\def\3n{\negthinspace \negthinspace \negthinspace }
\def\2n{\negthinspace \negthinspace }
\def\1n{\negthinspace }
\def\ns{\noalign{\smallskip} }

\def\ds{\displaystyle}
\def\G{\Gamma}
\def\D{\Delta}

\def\O{\Omega}
\def\cA{{\cal A}}

\def\cD{{\cal D}}
\def\cE{{\cal E}}
\def\cF{{\cal F}}
\def\cG{{\cal G}}

\def\cL{{\cal L}}

\def\cU{{\cal U}}

\def\cX{{\cal X}}

\def\cl{{\cal l}}
\def\no{\noindent}

\def\bs{\bigskip}
\def\q{\quad}
\def\qq{\qquad}

\def\liminf{\mathop{\underline{\rm lim}}}

\def\pa{\partial}

\def\wt{\widetilde}
\def\cd{\cdot}
\def\cds{\cdots}

\def\span{\hbox{\rm span$\,$}}

\def\cl{\overline}

\def\({\Big (}
\def\){\Big )}
\def\[{\Big[}
\def\]{\Big]}

\def\={\buildrel \triangle \over =}

\def\resp{{\it resp. }}

\def\be{\begin{equation}}
\def\bel{\begin{equation}\label}
\def\ee{\end{equation}}
\def\bea{\begin{eqnarray}}
\def\eea{\end{eqnarray}}
\def\bt{\begin{theorem}}
\def\et{\end{theorem}}
\def\bc{\begin{corollary}}
\def\ec{\end{corollary}}
\def\bl{\begin{lemma}}
\def\el{\end{lemma}}
\def\bp{\begin{proposition}}
\def\ep{\end{proposition}}
\def\br{\begin{remark}}
\def\er{\end{remark}}
\def\ba{\begin{array}}
\def\ea{\end{array}}
\def\bd{\begin{definition}}
\def\ed{\end{definition}}

\newtheorem{lemma}{Lemma}[section]
\newtheorem{remark}{Remark}[section]

\newtheorem{theorem}{Theorem}[section]
\newtheorem{corollary}{Corollary}[section]

\newtheorem{definition}{Definition}[section]
\newtheorem{proposition}{Proposition}[section]

\makeatletter
   
   \@addtoreset{equation}{section}
\makeatother

\begin{document}

\title{\bf  Averaged controllability for Random Evolution Partial Differential Equations}

\author{Qi L\"u\thanks{School of Mathematics,
Sichuan University, Chengdu, 610064, China
(\tt{luqi59@163.com}).} \q and \q Enrique
Zuazua\thanks{BCAM - Basque Center for Applied
Mathematics, Mazarredo, 14. E-48009
Bilbao-Basque
Country-Spain.}~\thanks{Ikerbasque, Basque
Foundation for Science, Mar'a Diaz de Haro 3,
E-48013 Bilbao-Basque Country-Spain
(\tt{zuazua@bcamath.org}).}}

\date{}

\maketitle

\begin{abstract}

We analyze the averaged controllability
properties of random evolution Partial
Differential Equations.

We mainly consider  heat and Schr\"odinger
equations with random parameters, although the
problem is also formulated in an abstract frame.

We show that the averages of parabolic equations
lead to parabolic-like dynamics that enjoy the
null-controllability properties of solutions of
heat equations in an arbitrarily short time and
from arbitrary measurable sets of positive
measure.

In the case of  Schr\"odinger equations we show
that, depending on the probability density
governing the random parameter, the average may
behave either as a conservative or a
parabolic-like evolution, leading to
controllability properties, in average, of very
different kind.
\end{abstract}

\bs


\bs

\no{\bf Key Words}. random evolution Partial
Differential Equations, averaged
controllability, averaged observability,
Schr\"odinger equation, heat equation.


\section{Introduction}


\q\;\,We analyze the problem of controlling
systems with randomly depending coefficients in
the context of evolution Partial Differential
Equations (PDEs). More precisely, we consider
the problem of averaged controllability which
consists, roughly, of controlling the averaged
dynamics, with respect to the random parameters.
This problem was introduced and solved in
\cite{Z4} in the context of finite dimensional
systems, where the same issue was also
formulated for  PDE.

When the dynamics of the state is governed by a
pair of random operators (determining the free
dynamics and the control operators,
respectively), generally speaking, controlling
the system would require to know  the actual
value of
 the random parameters. But this is unfeasible. To avoid this paradox,  controls should be independent of the
unknown parameters. But, of course, this
restricts our ability to deal with the
randomness of the system. Accordingly the
control requirement needs to be relaxed. The
most natural relaxation is to require the
control to perform optimally in an averaged
sense. This amounts to controlling the
mathematical expectation of the solutions,
 making a robust compromise of all the
possible realizations of the system for the
various possible values of the random
parameters.

In this work we consider mainly the heat and
Schr\"odinger equations, with diffusivity and
dispersivity operators depending on a random
variable in a multiplicative manner. We show
that, while the average of heat equations leads
also to a heat-like dynamics, the behavior of
the averages for the Schr\"odinger equations
depends in a very sensitive manner on the
density of probability of the random variable so
that, in some cases,  by averaging, this leads
to a dynamics of conservative nature, similar to
the original Schr\"odinger equation under
consideration and, in others, to a
parabolic-like behavior.

Our method of proof combines using Fourier
decomposition methods to identify the averaged
dynamics to, later,  utilizing the existing
tools developed  for the controllability of
parabolic and conservative systems, to deduce
the averaged controllability results.  When the
resulting averaged dynamics is of parabolic
nature, averaged null controllability is proved
for arbitrarily short time intervals and from
measurable sets of positive measure. On the
contrary,   when the averaged dynamics is of
conservative type, averaged exact
controllability is proved under suitable
geometric conditions on the support of the
controls that are by now well known in the
context of wave-like and Schr\"odinger-like
equations.

 In order to illustrate the effect of averaging and how it may change the dynamics of the original system,   let us consider the simplest transport equation
\begin{equation}\label{trans-system0}
\left\{
\begin{array}{ll}\ds
y_t + \a \cdot \nabla y = 0 &\mbox{ in }\dbR^d\times [0,\infty),\\
\ns\ds y(0)=y_0 &\mbox{ in }\dbR^d.
\end{array}
\right.
\end{equation}
Here $y_0\in L^2(\dbR^d)$ and
$\a(\cd):\O\to\dbR^d$ is a $d$-dimensional
standard normally distributed random variable,
with the probability density
$$
\rho(\a)=\frac{1}{(2\pi)^{\frac{d}{2}}}e^{-\frac{|\a|^2}{2}}
\mbox{ for } \a\in\dbR^d.
$$

The solution to \eqref{trans-system0} reads
$$
y(x,t,\o;y_0)=y_0(x-t\a) \mbox{ for
}(x,t)\in\dbR^d\times [0,\infty).
$$
Then, the mathematical expectation or averaged
state
$$
\begin{array}{ll}\ds
\tilde y(x,t)\3n&\ds\=\int_\O
y(x,t,\o;y_0)d\dbP(\o) =
\frac{1}{(2\pi)^{\frac{d}{2}}}\int_{\dbR^d}
y_0(x-\a t)e^{-\frac{|\a|^2}{2}}d\a
\\
\ns&\ds =
\frac{1}{(2\pi)^{\frac{d}{2}}t^d}\int_{\dbR^d}
y_0(z)e^{-\frac{(x-z)^2}{2t^2}}dz
\end{array}
$$
solves the following heat equation
\begin{equation}\label{6.14-eq2}
\left\{
\begin{array}{ll}\ds
\tilde y_t - \frac{1}{t}\D\tilde y= 0 &\mbox{ in
}
\dbR^d\times [0,\infty),\\
\ns\ds \tilde y(0) = y_0 &\mbox{ in }\dbR^d,
\end{array}
\right.
\end{equation}
namely, $\varphi(x, t)= \tilde y(x,\sqrt{2t})$
solves
\begin{equation}\label{6.29-eq1}
\left\{
\begin{array}{ll}\ds
\f_t -  \D\f = 0 &\mbox{ in }
\dbR^d\times [0,\infty),\\
\ns\ds \f(0) = y_0 &\mbox{ in }\dbR^d.
\end{array}
\right.
\end{equation}

The fact that averages of transport-like
equations may enjoy enhanced regularity
properties, first discovered in \cite{Ag}, is
well known in different contexts. There are some
different presentations of this smoothing
properties, referred to as ``averaging lemmas"
in the context of kinetic equations (see
\cite{BD,GLPS,GPS} for example). These smoothing
properties are useful when proving existence of
solutions of linear and nonlinear kinetic
equations \cite{BGPS,DP,DP1,LPT}. The proofs of
these ``averaging lemma" usually employ the
Fourier transform in all space-velocity-time
variables. In \cite{BD}, the authors gave a
proof that only uses the Fourier transform with
respect to $x$ and $v$, thus leading to a
Fourier representation of the average, evolving
in time,  which is very similar  in spirit to
the method we employ in this paper. Note however
that, for control purposes, we need to identify
the nature of the evolution associated with the
averaged variable, and not only its smoothing
properties.

\begin{remark}
The above computation shows that one can get
diffusion processes by averaging a simple random
convection process with respect to its velocity.
This is also well known from a different
perspective, in the context of chaotic and stiff
oscillatory systems, that can be regarded as the
characteristic systems of transport equations
(see \cite{CSSW,GKS,MS}).
\end{remark}

This example shows that, by averaging, the
solution of a random transport equation may lead
to a solution of a heat-like equation and,
consequently,  that time-reversible systems may
become strongly irreversible through averaging.
Furthermore, this occurs with the normally
distributed random variable which is ubiquitous
in nature due to the central limit theorem,
which states that the mean of many independent
random variables  drawn from the same
distribution is distributed approximately
normally, irrespective of the form of the
original distribution. Accordingly, in the real
world, physical quantities that are expected to
be the sum of many independent variables, such
as measurement errors, often have a distribution
very close to the normal distribution (see
 \cite{Br}).

The ``smoothing by averaging" effect mentioned
above has important consequences from a control
theoretical point of view as well.  Indeed,
while the transport equation
(\ref{trans-system0}), for a given value of the
random variable $\alpha$ (which determines the
velocity of propagation of waves),  enjoys the
property of exact controllability in finite
time, proportional to the travel time of
characteristics to get to the control set (the
boundary  or an open subset of the domain), the
averaged heat dynamics is controllable to zero
(or any other sufficiently smooth target) in an
arbitrarily short time and from any subset of
the domain where the dynamics evolves, without
any geometric condition on the support of the
controls, involving the propagation of
characteristics. Accordingly, through averaging,
we encounter on a single model, with randomly
depending coefficients, the classical dichotomy
arising in the context of controllability of
hyperbolic versus parabolic systems (see
\cite{Z2}).

This paper is devoted to systematically
addressing these questions in the context of
heat and Schr\"odinger equations. Our aim here
is not, by any means, to systematically address
all the possible scenarios but simply to
highlight some of the most fundamental phenomena
illustrating how, the existing tools for the
analysis of the controllability of PDE, can be
employed in this averaged context too. It is
important to highlight, however, that the
averaged states do not obey a PDE, not even a
semigroup. The dynamics can however be
represented in such a way that its main
controllability properties can be identified, by
analogy, with some of the main well-known
models, and analyzed by similar techniques.

In particular, we shall  show that the averages
of heat equations lead to heat-like dynamics
that enjoy the null-controllability properties
of solutions of heat equations in an arbitrarily
short time and from arbitrary measurable sets of
positive measure. In order to prove these
results we employ classical techniques based on
Carleman inequalities and the Fourier expansion
of solutions on the basis of the eigenfunctions
of the Laplacian generating the dynamics.

In the case of  Schr\"odinger equations we show
that, depending on the probability density, the
average may behave either as a conservative  or
a heat-like evolution, leading to
controllability properties, in average, of very
different kind. When the obtained average is of
parabolic nature the techniques above, employed
to treat the control properties of parabolic
averages, can be applied. However, when the
average behaves rather in a  conservative way we
employ specific techniques for the control of
wave-like equations.

The paper is organized as follows.  In section
\ref{abstract} we present all these problems in
an abstract setting in which different relevant
PDE models enter naturally.
 In Section \ref{sec-rheat}, we
study the null and approximate averaged
controllability problems for a class of random
heat equations. In Section \ref{sec-rsch}, we
study the null and exact averaged
controllability problem for a class of random
Schr\"odinger equations. In Section
\ref{sec-com} we give some further comments and
open problems.

\section{An abstract setting}\label{abstract}

\q\;\,Let $T>0$ and $E\subset [0,T]$ be a
Lebesgue measurable set with positive Lebesgue
measure. Let $H$ and $U$ be two Hilbert spaces.
Let $V\subset H$ be a Hilbert space which is
dense in $H$. Denote by $V'$ the dual space of
$V$ with respect to the pivot space $H$. Let
$(\Omega,\cF,\dbP)$ be a probability space. Let
$\{A(\o)\}_{\o\in\O}$ be a family of linear
operators satisfying the following conditions:

\begin{enumerate}
  \item $A(\cd)\in L^2(\O;\cL(D(A),H))$;
  \item $A(\o):D(A)\to H$ generates a $C_0$-semigroup
$\{S(t,\o)\}_{t\geq 0}$ on both $H$ and $V$ for
all $\o\in\O$;
  \item $S(t,\cd)y\in L^1(\O;V)$ for all
  $y\in V$ and $t\in [0,T]$.
\end{enumerate}
Let $B(\cd)\in L^2(\O;\cL(U,V))$.

Consider the following linear control system
\begin{equation}\label{ab-system0}
\left\{\begin{array}{ll} y_t(t)=A(\o)y(t)+ \chi_E(t)B(\o)u(t) & \mbox{ in } (0,T],\\
y(0)=y_0,
\end{array}\right.
\end{equation}
where $y_0\in V$ and  $u(\cd)\in L^2(E;U)$ is
the control.

In what follows, we denote by $y(\cd,\o;y_0)$
the solution to \eqref{ab-system0}, which is the
state of the system.  Although the initial datum
$y_0 \in V$ and the control $u(\cd)$ are
independent of the sample point $\o$, the state
$y(t, \o;y_0)$ of the system  depends on $\o$
nonlinearly.

According to the setting above,  for a.e.
$\o\in\O$, there is a solution $y(\cd,\o;y_0)\in
C([0,T];V)$ and the expectation or averaged
state $ \int_\O y(\cd,\o;y_0)d\dbP(\o) \in
C([0,T];V). $

We introduce the following notions of averaged
controllability for the system
\eqref{ab-system0}:

\begin{definition}\label{def-ex-con}
System \eqref{ab-system0} is said to fulfill the
property of {\it exact averaged controllability}
or to be  {\it exactly
 controllable in average} in  $E$ with control
cost $C>0$ if given any  $y_0, y_1\in V$, there
exists a control $u(\cd)\in L^2(E;U)$ such that
\begin{equation}\label{def-ex-con-eq1}
|u|_{L^2(E;U)}\leq C(|y_0|_V + |y_1|_V)
\end{equation}
and the average of solutions to
\eqref{ab-system0} satisfies
\begin{equation}\label{def-ex-con-eq2}
\int_\O y(T,\o;y_0) d\dbP(\o) =y_1.
\end{equation}
\end{definition}
\begin{remark}
The notion of exact averaged controllability was
first introduced in \cite{Z4}. A full
characterization was also given in the
finite-dimensional setting. In \cite{LZ},  this
issue was discussed for systems involving
finitely-many  linear parametric wave equations.
\end{remark}
\begin{definition}\label{def-nu-con}
System \eqref{ab-system0} fulfills the property
of {\it null averaged controllability} or  is
{\it null
 controllable in average} in $E$ with control
cost $C$ if given any initial datum $y_0 \in V$,
there exists a control $u\in L^2(E;U)$ such that
\begin{equation}\label{def-nu-con-eq1}
|u|_{L^2(E;U)}\leq C |y_0|_V
\end{equation}
and the average of the solutions to
\eqref{ab-system0} satisfies
\begin{equation}\label{def-nu-con-eq2}
\int_\O y(T,\o;y_0) d\dbP(\o) =0.
\end{equation}
\end{definition}
\begin{definition}\label{def-ap-con}
System \eqref{ab-system0} fulfills the property
of {\it approximate averaged controllability} or
is {\it approximately controllable in average}
in $E$ if given any $y_0, y_1 \in V$ and $\e>0$,
there exists a control $u_\e\in L^2(E;U)$ such
that the average of solutions to
\eqref{ab-system0} satisfies
$$
\Big|\int_\O y(T,\o;y_0)
d\dbP(\o)-y_1\Big|_V<\e.
$$
\end{definition}
\begin{remark}
As  in the finite dimensional context
(\cite{Z4}), we can also consider the averaged
control problem   with  random initial data,
i.e., $y_0\in L^2(\O;V)$. Nevertheless,
according to Remark \ref{rmk2}, this does not
lead to any essential new difficulty. Thus, for
the sake of simplicity of the presentation, we
only deal with the case where $y_0$ is
independent of $\o$.
\end{remark}

These notions are motivated by the problem of
controlling   a dynamics governed by a pair of
random operators $(A(\o),B(\o))$, where the
effective value of the parameter $\o$ is
unknown. Then, one aims at choosing a control,
independent of the unknown $\o$, to act
optimally in an averaged sense,
 making a robust compromise of all the
possible realizations of the system for the
various possible values of the sample point
$\o$. Similar problems can be considered in the
case where the initial datum to be controlled
depends on $\o$ too.

We have introduced the notions of
exact/null/approximate averaged controllability
in the framework of random evolution equations
but similar concepts make sense for parametrized
evolution equations (see \cite{Z4} for example).
In that context it is sufficient to replace the
expectation by a weighted average of the
parameter-depending controlled states.

\begin{remark} In the present context of randomly depending operators, the classical subordination properties of some control properties with respect to the others, that are classical for a given system,  have to be addressed more carefully. Of course, averaged exact controllability implies the averaged null and approximate controllability properties as well.
But, when $A$ and $B$ are independent of $\o$
and $A$ generates a $C_0$-group, exact
controllability is also a consequence of
 null controllability. However, the later
may  fail when considering averaged
controllability properties, as shown in the
example  in Remark \ref{rmk1} below.
\end{remark}
\begin{remark}\label{5.28-rem1}
For parametric control systems one can also
consider the problems of simultaneous
controllability and ensemble controllability,
which concern the possibility of controlling all
states with respect to different parameters
simultaneously by one single control. We refer
the readers to \cite{Li1,Li2} and \cite{BCP,LK}
for an introduction to these notions,
respectively. Of course the property of averaged
controllability we consider here is weaker than
these other ones since we only deal with the
average of the states with respect to the
parameters. But, as we shall see, averaged
controllability properties may be achieved  in
situations where simultaneous and ensemble
controllability are impossible. For instance,
let us consider the following control system:
\begin{equation}\label{5.28-eq16}
\left\{
\begin{array}{ll} y_t(t)=Ay(t)+  B(\o)u(t) & \mbox{ in } (0,T],\\
y(0)=y_0 \in \dbR^2,
\end{array}
\right.
\end{equation}
where
$$
A=\left(
    \begin{array}{cc}
      0 & 1 \\
      0 & 1 \\
    \end{array}
  \right), \qq B(\o) = B \mbox{ or } 2B \mbox{ for } B= \left(
                         \begin{array}{c}
                           0 \\
                           1 \\
                         \end{array}
                       \right).
$$
By Theorem 1 in \cite{Z4}, we know that the
system \eqref{5.28-eq16} is  null averaged
controllable.

But  system \eqref{5.28-eq16} is not
simultaneous null controllable.  Otherwise there
would exist a $u\in L^2(0,T)$ such that
$$
e^{AT}y_0 + \int_0^T e^{A(T-s)}Bu(s)ds =
e^{AT}y_0 + 2\int_0^T e^{A(T-s)}Bu(s)ds=0,
$$
which implies that $y_0=0$.
\end{remark}
\begin{remark}
The connection between averaged and simultaneous
controllability was analyzed in  \cite{LZ}
through  the following optimal control problem:

\medskip

Minimize $$J_k(u)\=\frac{1}{2}|u|_{L^2(0,T;U)}^2
+ k\int_\O|y(T)-y_1|^2 d\dbP(\o)$$ for all
$$
u\in \Big\{u\in L^2(0,T;U): \mbox{ The
corresponding solution $y$ satisfies that }
\int_\O y(T,\o)d\dbP(\o)=y_1\Big\}.
$$

\medskip

In that paper it is  proved that for every $k$,
$J_k(\cd)$ has a unique minimizer $u_k(\cd)$ and
that, if $H$ and $U$ are finite dimensional and
the system \eqref{ab-system0} is simultaneously
controllable, then $\{u_k\}_{k=1}^\infty$ weakly
converges to a simultaneous control.
\end{remark}
\begin{remark}
In the particular case that $A$ is independent
of $\o$,  the averaged controllability problems
can be reduced to the classical controllability
ones by setting
$$
\bar B = \int_\O B(\o) d\dbP(\o),\q \bar y(t) =
\int_\O y(t,\o;y_0) d\dbP(\o).
$$
Then we have that
\begin{equation}\label{6.8-eq4}
\left\{\begin{array}{ll} \bar y_t(t)=A\bar y(t)+\chi_E(t)\bar B u(t) &\mbox{ in } (0,T],\\
\bar y(0)=y_0.
\end{array}\right.
\end{equation}
The exact (\resp null,  approximate) averaged
controllability problems of \eqref{ab-system0}
are equivalent to the exact (\resp null,
approximate) controllability problem  of
\eqref{6.8-eq4}.
\end{remark}
\begin{remark}
Random evolution equations can be used to model
lots of uncertain physical processes (see
\cite{RSY,So,Sr} for example). Several notions
of controllability have been introduced but, as
far as we know, all of them concern driving the
state to a given destination by a control
depending on $\o$ (see \cite{JN,MR,MR1} and the
references therein). The property of averaged
controllability is, however, independent of the
specific realization of $\o$.
\end{remark}
Following the classical approach to deal with
controllability problems,  we introduce  the
adjoint system, which also depends on the
parameter $\o$:
\begin{equation}\label{ab-system1}
\left\{
\begin{array}{lll}
-z_t(t)=A^\ast (\o) z(t) & \mbox{ in } [0,T),\\
z(T)=z_0,
\end{array}
\right.
\end{equation}
where $z_0\in V'$.

Note that, in this adjoint system the initial
value (at time $t=T$) is taken to be independent
of $\o$. This is due to the fact that, although
$y(T,\o;y_0)$ depends on $\o$, its average,
which belongs to $V$, is of course independent
of $\o$. Then, to deal with the averaged state
it is sufficient to use as test functions,
adjoint states departing from configurations
that are independent of $\o$. This is the reason
we choose the final datum of \eqref{ab-system1}
to be independent of $\o$.

  As dual notions of the properties of averaged controllability  above we
introduce the following three concepts of
averaged observability.
\begin{definition}\label{def-ex-ob}
System \eqref{ab-system1} is exactly averaged
observable  or exactly observable in average in
$E$ if there is a constant $C>0$ such that for
any $z_0\in V'$,
\begin{equation}\label{def-ex-ob-eq1}
|z_0|^2_{V'} \leq C\int_0^T\chi_E(t)
\Big|\int_\O
B(\o)^*z(t,\o;z_0)d\dbP(\o)\Big|_{U}^2dt.
\end{equation}
\end{definition}
\begin{definition}\label{def-nu-ob}
System \eqref{ab-system1} is null averaged
observable or null observable in average in $E$
if there is a constant $C>0$ such that for any
$z_0\in H$,
\begin{equation}\label{def-nu-ob-eq1}
\Big|\int_\O z(0,\o;z_0)d\dbP(\o)\Big|_{V'}^2
\leq C\int_0^T \chi_E(t)\Big|\int_\O
B(\o)^*z(t,\o;z_0)d\dbP(\o)\Big|_{U}^2dt.
\end{equation}
\end{definition}
\begin{definition}\label{def-ap-ob}
System \eqref{ab-system1} is said to satisfy the
averaged unique continuation property in $E$ if
the fact that $\chi_E\int_\O
B(\o)^*z(\cd,\o;z_0)d\dbP(\o)=0$ in $L^2(0,T;U)$
implies that $z_0=0$.
\end{definition}

The adjoint system \eqref{ab-system1} being
exactly (\resp null) averaged observable means
that one can estimate the norm of the average of
the final (\resp initial) data of the adjoint
system, out of partial measurements done on the
averages of the adjoint states, with respect to
$\o$. These concepts have their own interest
when dealing with the observation of random
systems. The actual realization of the system
depending on $\o$ being unknown, it is natural
to address the problem based on the measurements
done on averages.

The weakest notion of averaged observability
under consideration is averaged unique
continuation. System \eqref{ab-system1}
satisfies the averaged unique continuation
property  when its state can be uniquely
determined by the partial measurements done on
the mathematical expectation. It is a natural
generalization of the unique continuation
property of evolution equations.

The average of the adjoint state, being
represented by $\int_\O z(t,\o;z_0) d\dbP(\o)$,
does not satisfy the semigroup property and it
is not a solution to an evolution equation.
Thus, one can not directly employ the existing
results on the observability of PDEs to
establish the averaged observability of the
adjoint system \eqref{def-ex-ob-eq1}. However,
as we shall see, by carefully analyzing and
identifying the dynamics generated by the
averages of the adjoint states, we shall be able
to apply the existing PDE techniques.

 In this paper, we mainly consider the
case that $A(\o)=\a(\o)A$ and $B(\o)=B$, where
$A$ generates a $C_0$-semigroup on $H$, $\a(\o)$
is a random variable and $B\in \cL(U,H)$. We
will show that the controllability properties of
the system \eqref{ab-system0} depend on the
choice of the random variable. We only consider
the following commonly used ones:

1. {\it Uniformly distributed random variable},
with probability density function $\rho(\cd)$ on
$[a,b]$, where $0<a<b$, $\a(\cd)$, given by
$$
\rho(\a)=\left\{
\begin{array}{ll}\ds\frac{1}{b-a}, &\mbox{ if
}\a\in [a,b],\\
\ns\ds 0, &\mbox{ if }\a\in (-\infty,a)\cup
(b,\infty).
\end{array}
\right.
$$

2.  {\it Exponentially distributed random
variable}, with probability density function
$\rho(\cd)$ given by
\begin{equation}\label{ed}
\rho(\a)= \left\{
\begin{array}{ll}\ds
e^{-(\a-c)}, &\mbox{ if } \a\geq c,\\
\ns\ds 0, &\mbox{ if } \a<c,
\end{array}
\right.
\end{equation}
for a given positive number $c$. In what
follows, for simplicity, we choose $c=1$ in
\eqref{ed}.

3.  {\it Standard normally distributed random
variable}, with probability density function of
$\rho(\cd)$ reads
$$
\rho(\a)=\frac{1}{\sqrt{2\pi}}e^{-\frac{\a^2}{2}}
\mbox{ for } \a\in\dbR.
$$

4. {\it  A random variable with standard Laplace
distribution}, with density function $\rho(\cd)$
given by
$$
\rho(\a)=\frac{1}{2}e^{-|\a|} \;\mbox{ for }
\;\a\in \dbR.
$$

5. {\it  A random variable with standard
Chi-squared distribution}, with density function
$\rho(\cd)$ given by
$$
\rho(\a;k)= \left\{
\begin{array}{ll}\ds
\frac{\a^{\frac{k}{2}-1}e^{-\frac{\a}{2}}}{2^{\frac{k}{2}}\G(\frac{k}{2})},
&\mbox{ if } \a\geq 0,\\
\ns\ds 0, &\mbox{ if } \a< 0,
\end{array}
\right.
$$
where $k\geq 1$ and
$\G(\frac{k}{2})=\int_0^\infty
\a^{\frac{k}{2}-1}e^{-\a}d\a$.

6. {\it  A random variable with standard Cauchy
distribution}, with density function $\rho(\cd)$
given by
$$
\rho(\a)=  \frac{1}{\pi(1+\a^2)}\q \mbox{ for
}\a\in\dbR.
$$
%


\section{Averaged controllability for random
 heat equations}\label{sec-rheat}


\q\;\,In this section, we study the null
averaged controllability problem for random heat
equations. Of course, the property of averaged
exact controllability has to be excluded because
of the very strong regularizing effect of heat
equations, that is preserved by averaging. Thus,
we focus on the properties of the  null and
approximate averaged controllability.

Let $T>0$  and $G\subset \dbR^d$ ($d\in\dbN$) be
a bounded domain with $C^2$ boundary $\pa G$ and
consider the following controlled random heat
equation:
\begin{equation}\label{heat-system0}
\left\{
\begin{array}{ll}\ds
y_t - \a\D y = \chi_{G_0\times E}u &\mbox{ in }
G\times
(0,T),\\
\ns\ds y=0 &\mbox{ on } \pa G\times
(0,T),\\
\ns\ds y(0)=y_0 &\mbox{ in }G.
\end{array}
\right.
\end{equation}
Here $y_0 \in L^2(G)$, and $G_0$ and $E$ are
subsets of $G$ and $[0, T]$ respectively, where
the controls are being applied.  The constant
diffusivity $\a:\O\to \dbR^+$ is assumed to be a
random variable. In this section, we make the
following assumptions on $G_0$ and $E$:

({\bf A1}) $G_0\subset G$ is a nonempty open
subset.

({\bf A2}) $E\subset [0,T]$ is a Lebesgue
measurable set with positive measure.

We first analyze the dynamics of the
mathematical expectation (average) of the
solutions of the random heat equations under
consideration. This will be done using Fourier
series expansions  and will give us an intuition
for the averaged observability and
controllability results to be expected and that
will be proved below.

Consider first  the following uncontrolled
random heat equation:
\begin{equation}\label{heat-system6}
\left\{
\begin{array}{ll}\ds
\hat y_t -  \a\D \hat y = 0 &\mbox{ in } G\times
(0,T),\\
\ns\ds \hat y=0 &\mbox{ on } \pa G\times
(0,T),\\
\ns\ds \hat y(0)=\hat y_0 &\mbox{ in }G.
\end{array}
\right.
\end{equation}
We decompose solutions in Fourier series on the
basis of the eigenfunctions of the Dirichlet
laplacian. To be more precise, consider the
unbounded linear operator $A_\D$ on
$L^2(G)$ given as%
$$
\left\{
\begin{array}{ll}\ds
D(A_\D) = H^2(G)\cap H_0^1(G),\\
\ns\ds A_\D u=-\D u, \q\mbox{ for any }u\in
D(A_\D).
\end{array}
\right.
$$
Let us denote  by $\{\l_j\}_{j=1}^\infty$ (with
$0<\l_1<\l_2\leq \cds$)  the eigenvalues of
$A_\D$ and let $\{e_j\}_{j=1}^\infty$ be the
corresponding eigenfunctions such that
$|e_j|_{L^2(G)}=1$ for $j\in\dbN$.

We assume that the initial datum takes the form
$\hat y_0=\sum_{j=1}^\infty \hat y_{0,j}e_j\in
L^2(G)$. The averaged state can also be
described in Fourier series as follows,
distinguishing the probability densities above.

{\bf Case 1}. If $\a(\cd)$ is a uniformly
distributed random variable on $[a,b]$, where
$0<a<b$, then,
\begin{equation}\label{avhe1}
\begin{array}{ll}\ds
\int_\O\hat y(x,t,\o;\hat y_0)d\dbP(\o)\3n&\ds=
\frac{1}{b-a}\int_a^b \sum_{i=1}^\infty \hat
y_{0,j}e^{-\l_j \a
t}e_jd\a\\
\ns&\ds = \frac{1}{b-a}\sum_{j=1}^\infty
\frac{1}{\l_j t}\hat y_{0,j}\big(e^{-\l_j a t} -
e^{-\l_j b t}\big)e_j.
\end{array}
\end{equation}
\begin{remark}
The values of a uniformly distributd random
variable, which is a relevant one in practice,
are uniformly distributed over an interval,
i.e., all points in the interval are equally
likely. It models the random phenomenon with
``equally possible outcomes".  When there is no
any a priori knowledge for $\a$ other than
$\a(\o)\in [a,b]$ for all $\o\in\O$, this is the
best possible choice. More details can be found
in \cite{JKB2}.
\end{remark}

{\bf Case 2.} When $\a(\cd)$ is the
exponentially distributed random variable,
\begin{equation}\label{aveh1}
\begin{array}{ll}\ds
\int_\O\hat y(x,t,\o;\hat y_0)d\dbP(\o)\3n&\ds=
\int_1^\infty e^{-(\a-1)}\sum_{j=1}^\infty \hat
y_{0,j}e^{-\l_j \a
t}e_jd\a\\
\ns&\ds = \sum_{j=1}^\infty \frac{1}{\l_j
t+1}\hat y_{0,j} e^{-\l_jt}e_j.
\end{array}
\end{equation}

\begin{remark}
The exponentially distributed random variable is
one of the most important random variables that
can be used to describe the time between events
in a process where they occur continuously and
independently at a constant average rate.
\end{remark}

In both  cases, the mathematical expectation of
the solution of the parameter-depending heat
equation evolves according to a heat-like
dynamics. The representation of the averages on
the basis of eigenfunctions exhibits the
exponential decay and smoothing effects that are
prototypical of heat-like problems.

Accordingly, the following null averaged
controllability result holds.

\begin{theorem}\label{con-th1}
Let ({\bf A1}) and ({\bf A2}) hold. Assume that,
either $\a(\cd)$ is a uniformly or exponentially
distributed random variable. Then, the system
\eqref{heat-system0} is null controllable in
average with control $u(\cd)\in
L^2(0,T;L^2(G_0))$. Further, there is a constant
$C>0$ such that
\begin{equation}\label{th1-eq1}
|u|_{L^2(0,T;L^2(G_0))} \leq C|y_0|_{L^2(G)}.
\end{equation}
\end{theorem}
\begin{remark}
In the context of heat equations, the null
controllability result with controls supported
in measurable sets is related to the bang-bang
property of the time optimal control  problems
with constrained controls (see \cite{PW,Wa} for
example).

More precisely, let us consider the following
heat equation:
\begin{equation}\label{8.18-eq1}
\left\{
\begin{array}{ll}\ds
\tilde y_t -\D \tilde y= \chi_{G_0}\tilde u
&\mbox{ in
} G\times [0,+\infty),\\
\ns\ds \tilde y=0 &\mbox{ on } \pa G\times
[0,+\infty),\\
\ns\ds \tilde y(0)=\tilde y_0 &\mbox{ in }G.
\end{array}
\right.
\end{equation}
Here the initial state $\tilde y_0\in L^2(G)$
and the control is assumed to belong to the
constrained set of admissible controls
$$
\tilde u\in \cU_M\=\{\tilde u\in
L^\infty(0,+\infty;L^2(G)):\,|\tilde
u|_{L^2(G)}\leq M \mbox{ a.e. }t\in
[0,+\infty)\}
$$
for some $M>0$. Let $\wt T^*$ be the
$\min\{\tilde t: \,\tilde y(\tilde t;\tilde
u,\tilde y_0)=0,\;\tilde u\in\cU_M\}$. A control
$u^*\in \cU_M$ such that $\tilde y(\wt
T^*;\tilde u,\tilde y_0)=0$ is called a time
optimal control and satisfies the bang-bang
property if $|\tilde u(t)|_{L^2(G)}=M$ for a.e.
$t\in [0,T^*]$. The proof of this fact requires
the null controllability with controls supported
in measurable sets (see \cite{Wa}).

The same bang-bang problem can be formulated in
the context of averaged null controllability we
are considering here. However, the techniques in
\cite{Wa} do not seem to apply because the
averages of the heat processes under
consideration do not satisfy the semigroup
property.

\end{remark}

To prove Theorem \ref{con-th1}, as in the
deterministic frame, we introduce the following
adjoint system:
\begin{equation}\label{heat-system1}
\left\{
\begin{array}{ll}\ds
z_t + \a\D z = 0 &\mbox{ in } G\times
(0,T),\\
\ns\ds z=0 &\mbox{ on } \pa G\times
(0,T),\\
\ns\ds z(T)=z_0 &\mbox{ in }G,
\end{array}
\right.
\end{equation}
where $z_0\in L^2(G)$. As mentioned above in the
abstract setting, the initial data (at time
$t=T$) of the adjoint system are assumed to be
independent of the random parameter. According
to Theorem \ref{th-nu-ob}, we only need to prove
that \eqref{heat-system1} is null observable in
average, which is a corollary of the following
result.

\begin{theorem}\label{th0}
Let ({\bf A1}) and  ({\bf A2}) hold. Assume that
$\a(\cd)$ is either an uniformly distributed or
an exponentially distributed random variable.
Then, there exists a constant $C>0$ such that
for any $y_0\in L^2(G)$ it holds that
\begin{equation}\label{th0-eq1}
\Big|\int_\O z(0,\o;y_0)d\dbP(\o)\Big|_{L^2(G)}
\leq C\int_E \Big|\int_\O
z(t,\o;y_0)d\dbP(\o)\Big|_{L^2(G_0)}dt.
\end{equation}
\end{theorem}

An immediate corollary of Theorem \ref{th0} is
as follows:
\begin{corollary}\label{cor1}
Let ({\bf A1}) and  ({\bf A2}) hold. Assume that
$\a(\cd)$ is either an uniformly distributed or
an exponentially distributed random variable.
Then the system \eqref{heat-system1} is null
observable in average.
\end{corollary}

{\it Proof}\,: From H\"older's inequality, we
have that
$$
\begin{array}{ll}\ds
\int_E \Big|\int_\O
z(t,\o;z_0)d\dbP(\o)\Big|_{L^2(G_0)}dt
\3n&\ds\leq \(\int_E dt\)^{\frac{1}{2}}\(\int_E
\Big|\int_\O
z(t,\o;z_0)d\dbP(\o)\Big|_{L^2(G_0)}^2dt\)^{\frac{1}{2}}\\
\ns&\ds \leq \sqrt{m(E)}\(\int_E \Big|\int_\O
z(t,\o;z_0)d\dbP(\o)\Big|_{L^2(G_0)}^2dt\)^{\frac{1}{2}}.
\end{array}
$$
This, together with the inequality
\eqref{th0-eq1}, implies that the system
\eqref{heat-system1} is null observable in
average.
\endpf

To prove Theorem \ref{th0}, we adopt the
strategy developed in \cite{Mi} in the context
of the heat equation, using spectral
decompositions. The null averaged observability
inequality is built in an iterative manner. More
precisely, we decompose the set $E$ into an
infinite sequence of connective (in time)
subsets in which an increasing number of Fourier
components of the average of the solution is
observed with uniform observability constants.
By iteration, the final datum is observed. To
apply this strategy, we need to use classical
results on how to divide $E$ into an infinite
sequence of subsets of positive Lebesgue
measure. We also need to know how to observe a
finite number of the Fourier components of the
average of the solution. These ingredients are
given in the following two lemmas.

\begin{lemma}\cite[Proposition 2.1]{PW}\label{lm1}
Let $E\subset[0,T]$ be a measurable set of
positive Lebesgue measure $m(E)$. Let $\ell$ be
a density point of $E$. Then for each $a
> 1$, there exists an $\ell_1 \in (\ell,T)$ such that
the sequence $\{\ell_k\}_{k=1}^\infty$, given by
\begin{equation}\label{lm1-eq1}
\ell_{k+1}=\ell+\frac{\ell_1-\ell}{a^k},
\end{equation}
satisfies
\begin{equation}\label{lm1-eq2}
m(E\cap (\ell_{k+1},\ell_k))\geq
\frac{\ell_k-\ell_{k+1}}{3}.
\end{equation}
\end{lemma}
\begin{lemma}\cite[Theorem 1.2]{Lu}\label{lm2}
There is a constant $C_1>0$ such that for any
$r>0$ and $\{a_j\}_{\l_j\leq r}\subset \dbC$,
\begin{equation}\label{lm2-eq1}
\(\sum_{\l_j\leq r} |a_j|^2\)^{\frac{1}{2}} \leq
C_1
e^{C_1\sqrt{r}}\(\int_{G_0}\Big|\sum_{\l_j\leq
r}a_j e_j(x) \Big|^2 dx\)^{\frac{1}{2}}.
\end{equation}
\end{lemma}
We are now in conditions to proceed to the proof
of Theorem \ref{th0}.

{\it Proof of Theorem \ref{th0}}\,: We only give
a proof for the case where $\a(\cd)$ is an
exponentially distributed random variable. The
proof for that $\a(\cd)$ is a uniformly
distributed random variable is very similar.

Put $\tilde z(x,t)=\int_\O
z(x,T-t,\o;z_0)d\dbP(\o)$. Let $z_{0,j} =
\langle z_0,e_j \rangle_{L^2(G)}$. Then, similar
to the computation of \eqref{aveh1}, we have
that
$$
\tilde z(x,t) =\sum_{j=1}^\infty \frac{1}{\l_j
t+1}z_{0,j} e^{-\l_jt}e_j.
$$
We only need to prove that
\begin{equation}\label{th1-eq2}
\Big|\sum_{j=1}^\infty z_{0,j}e^{-\l_j
T}e_j\Big|_{L^2(G)} \leq C\int_E |\tilde
z(x,t,\o) |_{L^2(G_0)}dt.
\end{equation}
Note that, in the right hand side of this
inequality, the observation is done in the
$L^1(E; L^2(G_0))$-norm. Thus, the result is
even stronger than the one we actually need, in
which the observation is done in $L^2(E;
L^2(G_0))$. As a consequence of the inequality
above we shall prove, the controls we shall
build will belong to $L^\infty(E; L^2(G_0))$.

Let $X_r=\span\{e_j\}_{\l_j\leq r}$ for each
$r>0$. For any $\xi\in L^2(G)$, we put
\begin{equation}\label{5.29-eq3}
S(t,\xi)=\sum_{j=1}^\infty \frac{1}{\l_j
t+1}\xi_{j} e^{-\l_jt}e_j,
\end{equation}
where $\xi_j=\langle \xi, e_j\rangle_{L^2(G)}$.
Then we see that
\begin{equation}\label{5.29-eq5}
\big|S(t,\xi)\big|_{L^2(G)}^2 \leq
\big|S(s,\xi)\big|_{L^2(G)}^2, \q\mbox{ for }
0\leq s\leq t\leq T,
\end{equation}
and
\begin{equation}\label{5.29-eq4}
\big|S(t,\xi)\big|_{L^2(G)}^2 \leq e^{-r
(t-s)}|S(s,\xi)|_{L^2(G)}, \q\mbox{ for all
}\xi\in X_r^\perp \mbox{ and } 0\leq s\leq t\leq
T.
\end{equation}
Let $\ell$ be a density point for $E$. By Lemma
\ref{lm1}, for a given $a > 1$, there exists   a
sequence $\{\ell_k\}_{k=1}^\infty$ satisfying
\eqref{lm1-eq1} and \eqref{lm1-eq2}.

We now define a sequence of subsets
$\{E_k\}_{k=1}^\infty$ of $(0,T)$ in the
following way:
\begin{equation}\label{5.28-eq1}
E_k\=\Big\{t-\frac{\ell_k-\ell_{k+1}}{6}:\, t\in
E\cap
\(\ell_{k+1}+\frac{\ell_k-\ell_{k+1}}{6},\ell_k\)\Big\},
\;\mbox{ for } k\in\dbN.
\end{equation}
Clearly, $E_k\subset (\ell_{k+1},\ell_{k+1} +
\frac{5}{6}(\ell_k-\ell_{k+1}))$. From
\eqref{lm1-eq2}, we have that
\begin{equation}\label{5.28-eq2}
\begin{array}{ll}\ds
m(E_k)\3n&\ds =m\(E\cap
\(\ell_{k+1}+\frac{\ell_k-\ell_{k+1}}{6},\ell_k\)\)
\\
\ns&\ds=m\(E\cap \[ (\ell_{k+1},\ell_k)\setminus
\(\ell_{k+1},\ell_{k+1}+\frac{\ell_k-\ell_{k+1}}{6}\)
\]\)\\
\ns&\ds \geq m(E\cap
(\ell_{k+1},\ell_k))-\frac{\ell_k-\ell_{k+1}}{6}
\\
\ns&\ds \geq \frac{\ell_k-\ell_{k+1}}{6}.
\end{array}
\end{equation}
Let $b > a $ be a positive number such that
\begin{equation}\label{5.28-eq3}
\frac{b}{a }> \frac{C_1 +
6\ln(12C_1a)}{(a-1)(\ell_1-\ell)}.
\end{equation}
Set $r_k = b^{2k}$. From \eqref{5.29-eq5}, we
have that for any $\xi\in X_{r_k}$,
$$
\begin{array}{ll}\ds
\q\int_{\ell_{k+1}}^{\ell_{k+1}+\frac{5}{6}(\ell_k-\ell_{k+1})}
\chi_{E_k}(t) \Big|
S\(\ell_{k+1}+\frac{5}{6}(\ell_k-\ell_{k+1}),\xi\)
\Big|_{L^2(G)}dt \\
\ns\ds\leq
\int_{\ell_{k+1}}^{\ell_{k+1}+\frac{5}{6}(\ell_k-\ell_{k+1})}\chi_{E_k}(t)
\big|S(t,\xi) \big|_{L^2(G)}dt.
\end{array}
$$
From \eqref{lm2-eq1} and \eqref{5.28-eq2}, we
find that
\begin{equation}\label{5.28-eq4}
\begin{array}{ll}\ds
\q\frac{\ell_k-\ell_{k+1}}{6}\Big| S
\(\ell_{k+1}+\frac{5}{6}(\ell_k-\ell_{k+1}),\xi\)
\Big|_{L^2(G)}\\
\ns\ds \leq m(E_k)\Big| S
\(\ell_{k+1}+\frac{5}{6}(\ell_k-\ell_{k+1}),\xi\)
\Big|_{L^2(G)}\\
\ns\ds \leq
\int_{\ell_{k+1}}^{\ell_{k+1}+\frac{5}{6}(\ell_k-\ell_{k+1})}\chi_{E_k}(t)
\big|S(t,\xi)\big|_{L^2(G)}dt\\
\ns\ds \leq
C_1e^{C_1\sqrt{r_k}}\int_{\ell_{k+1}}^{\ell_{k+1}+\frac{5}{6}(\ell_k-\ell_{k+1})}\chi_{E_k}(t)
\big|S(t,\xi)\big|_{L^2(G_0)}dt.
\end{array}
\end{equation}
Let $z_0 = z_0^1 + z_0^2$, where $z_0^1\in
X_{r_k}$ and $z_0^2\in X_{r_k}^\perp$. Taking
$\xi = S\(\frac{\ell_k-\ell_{k+1}}{6},z_0^1\)$
in \eqref{5.28-eq4}, we get that
\begin{equation}\label{5.28-eq5}
\begin{array}{ll}\ds
\q\frac{\ell_k-\ell_{k+1}}{6}|S(\ell_k,z_0^1)|_{L^2(G)}\\
\ns\ds \leq
\int_{\ell_{k+1}}^{\ell_{k+1}+\frac{5}{6}(\ell_k-\ell_{k+1})}
\chi_{E_k}(t)\Big|S\(t+\frac{\ell_k-\ell_{k+1}}{6},z_0^1\)\Big|_{L^2(G)}dt\\
\ns\ds \leq C_1
e^{C_1\sqrt{r_k}}\int_{\ell_{k+1}}^{\ell_{k+1}+\frac{5}{6}(\ell_k-\ell_{k+1})}
\chi_{E_k}(t)\Big|S\(t+\frac{\ell_k-\ell_{k+1}}{6},z_0^1\)\Big|_{L^2(G_0)}dt\\
\ns\ds \leq C_1
e^{C_1\sqrt{r_k}}\int_{\ell_{k+1}+\frac{\ell_k-\ell_{k+1}}{6}}^{\ell_k}\chi_{E_k}\(t-\frac{\ell_k-\ell_{k+1}}{6}\)
|S(t,z_0^1)|_{L^2(G_0)}dt.
\end{array}
\end{equation}
By \eqref{5.28-eq1}, we have that
\begin{equation}\label{5.28-eq13}
\chi_{E_k}\(t-\frac{\ell_k-\ell_{k+1}}{6}\)=\chi_E(t),
\q\mbox{ for any }t\in
\(\ell_{k+1}+\frac{\ell_k-\ell_{k+1}}{6},\ell_k\).
\end{equation}
Combining \eqref{5.29-eq5}, \eqref{5.29-eq4},
\eqref{5.28-eq5} and \eqref{5.28-eq13}, we find
that
\begin{equation}\label{5.28-eq6}
\begin{array}{ll}\ds
\q\frac{\ell_k-\ell_{k+1}}{6}|S(\ell_k,z_0^1)|_{L^2(G)}\\
\ns\ds \leq C_1
e^{C_1\sqrt{r_k}}\int_{\ell_{k+1}+\frac{\ell_k-\ell_{k+1}}{6}}^{\ell_k}
\chi_E(t) |S(t,z_0^1)|_{L^2(G_0)}dt
\\
\ns\ds \leq C_1
e^{C_1\sqrt{r_k}}\int_{\ell_{k+1}+\frac{\ell_k-\ell_{k+1}}{6}}^{\ell_k}
\chi_E(t)\big(|S(t,z_0)|_{L^2(G_0)} +
|S(t,z_0^2)|_{L^2(G)}\big)dt\\
\ns\ds \leq C_1
e^{C_1\sqrt{r_k}}\int_{\ell_{k+1}+\frac{\ell_k-\ell_{k+1}}{6}}^{\ell_k}
\chi_E(t) |S(t,z_0)|_{L^2(G_0)} dt\\
\ns\ds \q + C_1
e^{C_1\sqrt{r_k}}(\ell_k-\ell_{k+1})\Big|
S\(\ell_{k+1}+\frac{\ell_k-\ell_{k+1}}{6},z_0^2\)
\Big|_{L^2(G)}  \\
\ns\ds \leq C_1
e^{C_1\sqrt{r_k}}\int_{\ell_{k+1}}^{\ell_k}
\chi_E(t) |S(t,z_0)|_{L^2(G_0)} dt \\
\ns\ds \q + C_1
e^{C_1\sqrt{r_k}}(\ell_k-\ell_{k+1})e^{-
\frac{\ell_k-\ell_{k+1}}{6}r_k}\big| S
(\ell_{k+1},z_0^2) \big|_{L^2(G)}.
\end{array}
\end{equation}
Therefore, we obtain that
\begin{equation}\label{5.28-eq8}
\begin{array}{ll}\ds
\q\frac{\ell_k-\ell_{k+1}}{6}|S(\ell_k,z_0)|_{L^2(G)}\\
\ns\ds \leq C_1
e^{C_1\!\sqrt{r_k}}\int_{\ell_{k+1}}^{\ell_k}\!\!
\chi_E(t) |S(t,z_0)|_{L^2(G_0)} dt\! +\! C_1
e^{C_1\!\sqrt{r_k}}(\ell_k\!-\!\ell_{k+1})e^{-
\frac{\ell_k-\ell_{k+1}}{6}r_k}\big|
S(\ell_{k+1},z_0^2)\big|_{L^2(G)}\\
\ns\ds \q +
\frac{\ell_k-\ell_{k+1}}{6}|S(\ell_k,z_0^2)|_{L^2(G)}\\
\ns\ds \leq C_1
e^{C_1\!\sqrt{r_k}}\int_{\ell_{k+1}}^{\ell_k}\!\!
\chi_E(t) |S(t,z_0)|_{L^2(G_0)} dt\! +\! C_1
e^{C_1\!\sqrt{r_k}}(\ell_k\!-\!\ell_{k+1})e^{-
\frac{\ell_k-\ell_{k+1}}{6}r_k}\big|
S(\ell_{k+1},z_0^2) \big|_{L^2(G)}\\
\ns\ds \q +
\frac{\ell_k-\ell_{k+1}}{6}e^{-r_{k}(\ell_k-\ell_{k+1})}|S(\ell_{k+1},z_0^2)|_{L^2(G)}.
\end{array}
\end{equation}
Thus, it holds that
\begin{equation}\label{5.28-eq9}
\begin{array}{ll}\ds
\q\frac{\ell_k-\ell_{k+1}}{6}|S(\ell_k,z_0)|_{L^2(G)}\\
\ns\ds \leq C_1
e^{C_1\sqrt{r_k}}\int_{\ell_{k+1}}^{\ell_k}
\chi_E(t) |S(t,z_0)|_{L^2(G_0)} dt\\
\ns\ds \q + (\ell_k-\ell_{k+1})e^{-
\frac{\ell_k-\ell_{k+1}}{6}r_k} \big(C_1
e^{C_1\sqrt{r_k}} + 1\big) \big|
S(\ell_{k+1},z_0)\big|_{L^2(G)}.
\end{array}
\end{equation}
This concludes that
\begin{equation}\label{5.28-eq10}
\begin{array}{ll}\ds
\frac{\ell_k-\ell_{k+1}}{6C_1e^{C_1\sqrt{r_k}}}|S(\ell_{k},z_0)|_{L^2(G)}
-
\frac{C_1e^{C_1\sqrt{r_k}}+1}{C_1e^{C_1\sqrt{r_k}}}
(\ell_k-\ell_{k+1})e^{-
\frac{\ell_k-\ell_{k+1}}{6}r_k}|S(\ell_{k+1},z_0)|_{L^2(G)}\\
\ns\ds \leq
\int_{\ell_{k+1}}^{\ell_k}\chi_E(t)|S(t,z_0)|_{L^2(G_0)}dt.
\end{array}
\end{equation}
By summing the inequality \eqref{5.28-eq10} from
$k=1$ to $k=\infty$, we obtain that
\begin{equation}\label{5.28-eq11}
\frac{\ell_{1}-\ell_{2}}{6C_1e^{C_1\sqrt{r_{1}}}}
|S(\ell_{1},z_0)|_{L^2(G)} + \sum_{k=1}^\infty
f_k |S(\ell_{k+1},z_0)|_{L^2(G)}\leq \int_0^T
\chi_E(t)|S(t,z_0)|_{L^2(G_0)}dt,
\end{equation}
where
$$
f_k=\frac{\ell_{k+1}-\ell_{k+2}}{6C_1e^{C_1\sqrt{r_{k+1}}}}-
\frac{C_1e^{C_1\sqrt{r_k}}+1}{C_1e^{C_1\sqrt{r_k}}}(\ell_k-\ell_{k+1})e^{-
\frac{\ell_k-\ell_{k+1}}{6}r_k}, \q k=1,2,\cds
$$
From  \eqref{5.28-eq3} and $r_k = b^{2k}$, we
have that
$$
f_k\geq 0 \;\mbox{ for any }\; k=1, 2, \cds
$$
This, together with \eqref{5.28-eq11}, deduces
that
\begin{equation}\label{5.28-eq12}
|S(\ell_{1},z_0)|_{L^2(G)}\leq
\frac{6C_1e^{C_1\sqrt{r_{1}}}}{\ell_{1}-\ell_{2}}\int_E|S(t,z_0)|_{L^2(G_0)}dt.
\end{equation}
Since $\ell_1<T$, we can find a constant $C>0$
such that
\begin{equation}\label{5.28-eq14}
\frac{C}{1+\l_j\ell_1}\geq e^{-\l_j(T-\ell_1)}
\mbox{ for every } j\in\dbN.
\end{equation}
From \eqref{5.28-eq12} and \eqref{5.28-eq14}, we
get that
\begin{equation}\label{5.28-eq15}
\Big|\sum_{j=1}^\infty z_{0,j}e^{-\l_j
T}e_j\Big|_{L^2(G)}\leq
C\int_E|S(t,z_0)|_{L^2(G_0)}dt.
\end{equation}
This completes the proof.
\endpf

 As an easy corollary of Theorem \ref{th0},
we have the following result.

\begin{theorem}\label{th-ap-con1}
Let ({\bf A1}) and ({\bf A2}) hold.  System
\eqref{heat-system0} is approximately
controllable in average, provided that $\a(\cd)$
is a uniformly distributed or an exponentially
distributed random variable.
\end{theorem}

{\it Proof}\,: According to Theorem
\ref{th-ap-ob}, we only need to prove that the
unique continuation property in $E$ is
satisfied.

Assume that $z=0$ in $G_0\times E$. From Theorem
\ref{th0},  we obtain that
$$
\Big|\int_\O z(0,\o;z_0)d\dbP(\o)\Big|_{L^2(G)}
\leq C\int_E \Big|\int_\O
z(t,\o;z_0)d\dbP(\o)\Big|_{L^2(G_0)}dt.
$$
Hence,  $\int_\O z(0,\o;z_0)d\dbP(\o)=0$.
 On the other hand, if $z_{0,j}=\int_G z_0 e_jdx$,
\begin{equation}\label{8.12-eq1}
0=\int_\O z(0,\o;z_0)d\dbP(\o)=\int_1^\infty
e^{-(\a-1)}\sum_{j=1}^\infty z_{0,j}e^{-\l_j \a
T}e_jd\a =\sum_{j=1}^\infty \frac{1}{\l_j
T+1}z_{0,j} e^{-\l_jT}e_j.
\end{equation}
Since $\{e_j\}_{j=1}^\infty$ is an orthonormal
basis of $L^2(G)$, it follows that
$\frac{1}{\l_j T+1}z_{0,j} e^{-\l_jT}=0$ for all
$j\in\dbN$. Thus, we get that $z_{0,j}=0$ for
all $j\in\dbN$, which implies that $z_0=0$.
\endpf


\section{Averaged controllability for random
Schr\"odinger equations}\label{sec-rsch}

\subsection{Preliminaries}

\q\;\,In this section, we study the null and
exact averaged controllability problems for a
class of random Schr\"odinger equations of the
form
\begin{equation}\label{sch-system0}
\left\{
\begin{array}{ll}\ds
y_t - i\a\D y = \chi_{G_0\times E}u &\mbox{ in }
G\times
(0,T),\\
\ns\ds y=0 &\mbox{ on } \pa G\times
(0,T),\\
\ns\ds y(0)=y_0 &\mbox{ in }G.
\end{array}
\right.
\end{equation}
Here $\a(\cd):\O\to\dbR$ is a random variable,
the initial datum $y_0$ belongs to $L^2(G)$ and
$G_0$ is a suitable subdomain of $G$.

In this time-reversible setting the
Schr\"odinger equation is well-posed whatever
the sign of $\a$ is, contrary to the   heat
equation. Thus,  we have more choices for the
random variable $\a(\cd)$.

The average of the solutions to random
Schr\"odinger equations may lead to very
different dynamics, depending on the random
variable under consideration. To see this, let
us first consider the Schr\"odinger system in
the absence of control:
\begin{equation}\label{sch-system5}
\left\{
\begin{array}{ll}\ds
\hat y_t - i\a\D \hat y = 0 &\mbox{ in } G\times
(0,T),\\
\ns\ds \hat y=0 &\mbox{ on } \pa G\times
(0,T),\\
\ns\ds \hat y(0)=\hat y_0 &\mbox{ in }G.
\end{array}
\right.
\end{equation}
Here $\hat y_0=\sum_{j=1}^\infty \hat
y_{0,j}e_j\in L^2(G)$.

{\bf Case 1}. When $\a(\cd)$ is a uniformly
distributed random variable on $[a,b]$, where
$a,b\in\dbR$, then,
\begin{equation}\label{avu1}
\begin{array}{ll}\ds
\int_\O\hat y(x,t,\o;\hat y_0)d\dbP(\o)\3n&\ds=
\frac{1}{b-a}\int_a^b \sum_{i=1}^\infty \hat
y_{0,j}e^{-i\l_j \a
t}e_jd\a\\
\ns&\ds = \frac{1}{b-a}\sum_{j=1}^\infty
\frac{1}{i\l_j t}\hat y_{0,j}\big(e^{-i\l_j a t}
- e^{-i\l_j b t}\big)e_j.
\end{array}
\end{equation}
{\bf Case 2}.  When $\a(\cd)$ is an
exponentially distributed random variable, then
\begin{equation}\label{ave1}
\begin{array}{ll}\ds
\int_\O\hat y(x,t,\o;\hat y_0)d\dbP(\o)\3n&\ds=
\int_1^\infty e^{-(\a-1)}\sum_{j=1}^\infty \hat
y_{0,j}e^{-i\l_j \a
t}e_jd\a\\
\ns&\ds = \sum_{j=1}^\infty \frac{1}{i\l_j
t+1}\hat y_{0,j} e^{-i\l_jt}e_j.
\end{array}
\end{equation}
{\bf Case 3}.  For the  normally distributed
random variable $\a(\cd)$ we have,
\begin{equation}\label{ave1.1}
\begin{array}{ll}\ds
\int_\O\hat y(x,t,\o;\hat y_0)d\dbP(\o)\3n&\ds =
\frac{1}{\sqrt{2\pi}}\int_{-\infty}^\infty
e^{-\frac{\a^2}{2}}\sum_{j=1}^\infty \hat
y_{0,j}e^{-i\l_j \a
t}e_jd\a\\
\ns&\ds= \sum_{j=1}^\infty \hat
y_{0,j}e^{-\frac{1}{2}\l_j^2 t^2}e_j.
\end{array}
\end{equation}
{\bf Case 4}.  When $\a(\cd)$ is a random
variable with Laplace distribution we have
\begin{equation}\label{ave1.2}
\begin{array}{ll}\ds
\int_\O\hat y(x,t,\o;\hat y_0)d\dbP(\o)\3n&\ds =
\frac{1}{2}\int_{-\infty}^\infty
e^{-|\a|}\sum_{j=1}^\infty \hat y_{0,j}e^{-i\l_j
\a
t}e_jd\a\\
\ns&\ds= \sum_{j=1}^\infty \frac{1}{1+\l_j^2
t^2}\hat y_{0,j}e^{-i\l_j t}e_j.
\end{array}
\end{equation}
\begin{remark}
The Laplace distribution can be thought of as
two exponential distributions (with an
additional location parameter) spliced together
back-to-back. It governs the difference of two
independent identically distributed exponential
random variables and can be regarded as the
generalization of the exponential distribution
on the whole real line. More details can be
found in \cite{KKP}.
\end{remark}

{\bf Case 5}.  For the Chi-squared distribution
it holds:
\begin{equation}\label{ave1.3}
\begin{array}{ll}\ds
\int_\O\hat y(x,t,\o;\hat y_0)d\dbP(\o)\3n&\ds =
\int_{0}^\infty
\frac{\a^{\frac{k}{2}-1}e^{-\frac{\a}{2}}}{2^{\frac{k}{2}}\G(\frac{k}{2})}\sum_{j=1}^\infty
\hat y_{0,j}e^{-i\l_j \a
t}e_jd\a\\
\ns&\ds= \sum_{j=1}^\infty \frac{1}{(1+2i\l_j
t)^{\frac{k}{2}}}\hat y_{0,j}e^{-i\l_j t}e_j.
\end{array}
\end{equation}
\begin{remark}
A Chi-squared distributed random variable is the
sum of the squares of $k$ independent standard
normally distributed random variables. It is one
of the most widely used probability
distributions in inferential statistics. We
refer the readers to \cite{JKB1} for more
details.
\end{remark}

{\bf Case 6}.  For the Cauchy distribution we
have
\begin{equation}\label{ave1.4}
\begin{array}{ll}\ds
\int_\O\hat y(x,t,\o;\hat y_0)d\dbP(\o)\3n&\ds =
\int_{-\infty}^\infty
\frac{1}{\pi(1+\a^2)}\sum_{j=1}^\infty \hat
y_{0,j}e^{-i\l_j \a
t}e_jd\a\\
\ns&\ds= \sum_{j=1}^\infty \hat y_{0,j}e^{-\l_j
t}e_j.
\end{array}
\end{equation}
\begin{remark}
The Cauchy distribution is associated with many
processes, including resonance energy
distribution, impact and natural spectral and
quadratic stark line broadening. It also has
important connections with other random
variables. For example, when $\g_1$ and $\g_2$
are two independent standard normally
distributed random variables, then the ratio
$\g_1/\g_2$ has the standard Cauchy
distribution. More details can be found in
\cite{JKB1}.
\end{remark}
\begin{remark}
From \eqref{ave1.4}, we know that when $\a$ is a
random variable with Cauchy distribution, then
the average of the solution to the random
Schr\"odinger equation \eqref{sch-system5}
becomes a solution of the heat equation. This is
another example that after averaging, one enjoy
enhanced regularity properties.
\end{remark}

According to the above results, there are
essentially two different dynamics for the
 averages, depending on whether the time-exponentials entering in the Fourier expansion are real or imaginary. In cases 3 and 6, the
average has a heat-like behavior. However, in
cases 1, 2, 4 and 5, the average has a
Schr\"odinger-like behavior.


\subsection{Null averaged controllability }


\q\; We first recall the following assumptions
on $G_0$ and $E$:

({\bf A1}) Let $G_0\subset G$ be a nonempty open
subset.

({\bf A2}) Let $E\subset [0,T]$ be a Lebesgue
measurable set with positive measure.

\begin{theorem}\label{con-th2}
Let ({\bf A1}) and ({\bf A2}) hold. If $\a$ is a
random variable with normal distribution or
Cauchy distribution, then the system
\eqref{sch-system0} is null controllable in
average with control $u(\cd)\in
L^2(0,T;L^2(G_0))$. Further, there is a constant
$C>0$ such that
\begin{equation}\label{cth2-eq1}
|u|_{L^2(0,T;L^2(G_0))} \leq C|y_0|_{L^2(G)}.
\end{equation}
\end{theorem}
\begin{remark}
Note that, in the present case,  the random
Schr\"odinger equations is null controllable in
average without any assumption on the support
$G_0$ of the control, other than being of
positive measure. This is in contrast with the
well known results on the null controllability
of Schr\"odinger equations, where $G_0$ is
assumed, for instance, to fulfill the classical
Geometric Control Condition (GCC)(see
\cite{Le1,Le} for example) or other geometric
restrictions associated to multiplier techniques
or Carleman inequalities(see \cite{LTZ,Ma} for
example). In the present case, these
restrictions on $G_0$ are not needed since the
averages behave in a parabolic fashion.

\end{remark}

To prove Theorem \ref{con-th2}, we introduce the
adjoint system of \eqref{sch-system0} as
follows:
\begin{equation}\label{sch-system1}
\left\{
\begin{array}{ll}\ds
z_t + i\a\D z = 0 &\mbox{ in } G\times
(0,T),\\
\ns\ds z=0 &\mbox{ on } \pa G\times
(0,T),\\
\ns\ds z(T)=z_0 &\mbox{ in }G.
\end{array}
\right.
\end{equation}
By Theorem \ref{th-nu-ob}, we only need to prove
the following result.
\begin{theorem}\label{ob-th1}
Under the assumptions of Theorem \ref{con-th2}
the system \eqref{sch-system1} is null
observable in average if $\a$ is a random
variable with normal or Cauchy distribution.
\end{theorem}

Indeed, we have the following stronger
observability estimates.

\begin{proposition}\label{th2}
Let ({\bf A1}) and ({\bf A2}) hold. Assume that
$\a(\cd)$ is a random variable with normal
distribution or Cauchy distribution. Then there
exists a constant $C>0$ such that for any
$z_0\in L^2(G)$, it holds that
\begin{equation}\label{th2-eq1}
\Big|\int_\O z(0,\o;z_0)d\dbP(\o)\Big|_{L^2(G)}
\leq C\int_E \Big|\int_\O
z(x,t,\o;z_0)d\dbP(\o)\Big|_{L^2(G_0)}dt.
\end{equation}
\end{proposition}

The proof of Proposition \ref{th2} is very
similar to the one for Theorem \ref{th0}. We
omit it here.

\begin{remark}\label{rmk1}

If $\a(\cd)$ is a random variable with Cauchy
distribution, then we have that
\begin{equation}\label{savc1}
\begin{array}{ll}\ds
\int_\O z(x,t,\o;z_0)d\dbP(\o)\3n&\ds =
\int_{-\infty}^\infty
\frac{1}{\pi(1+\a^2)}\sum_{j=1}^\infty
z_{0,j}e^{-i\l_j \a
(T-t)}e_jd\a\\
\ns&\ds= \sum_{j=1}^\infty z_{0,j}e^{-\l_j
(T-t)}e_j.
\end{array}
\end{equation}
In this case, Proposition \ref{th2} is an
immediate corollary of the observability
estimate for heat equations.

This is an example of a system that is null but
not exactly controllable in average.
\end{remark}
%


\subsection{Exact averaged controllability}


In this subsection we consider the cases where
the averages behave as Schr\"odinger-like
semgroups.  We assume that $E=[0,T]$ so that the
control is active in any time instant within the
time interval $[0, T]$.

Let us first consider the following
Schr\"odinger equation
\begin{equation}\label{Sch-system2}
\left\{
\begin{array}{ll}\ds
\f_t + \k i\D \f = 0 &\mbox{ in } G\times
(0,T],\\
\ns\ds \f = 0 &\mbox{ on } \pa G\times
(0,T),\\
\ns\ds \f(0)=\f_0 &\mbox{ in }G,
\end{array}
\right.
\end{equation}
where $\k \in \dbR\setminus\{0\}$ and $\f_0\in
L^2(G)$.

We make the following assumption in this
subsection on $G_0$:

({\bf A3}) Whatever $T>0$ and $k \ne 0$ are,
there is a constant $C>0$ such that for any
$\f_0\in L^2(G)$, the solution $\f(\cd,\cd)$ to
\eqref{Sch-system2} satisfies
\begin{equation}\label{6.29-eq2}
|\f_0|_{L^2(G)}^2\leq C\int_0^T\int_{G_0}
|\f|^2dxdt.
\end{equation}
We refer the readers to \cite{An,Ja,Le1,Ph} for
the study of \eqref{6.29-eq2} under different
conditions for $G_0$. In those articles one can
find various sufficient conditions on the subset
$G_0$ so that the observability inequality above
holds, depending of the techniques of proof
employed (multipliers, Carleman inequalities,
Microlocal analysis). In particular, this
observability inequality for the Schr\"odinger
equation holds as soon as it is satisfied for
the wave equation in some time horizon. Thus, in
particular, it holds under the classical
Geometric Control Condition (GCC) guaranteeing,
roughly, that all rays of geometric optics enter
the observation set $G_0$ in some uniform time.

We have the following observability result for
the system \eqref{Sch-system2}.

\begin{theorem}\label{con-th6}
The following results hold:
\begin{itemize}
  \item Let $\a(\cd)$ be a uniformly distributed
random variable on an interval $[a,b]$. Then,
the system \eqref{sch-system0} is exactly
controllable in average with $V=H=L^2(G)$ and
$U=H^{-2}(G_0)$.

  \item Let $\a(\cd)$ be an exponentially
distributed random variable. Then, the system
\eqref{sch-system0} is exactly controllable in
average with $H=L^2(G)$, $V=H^2(G)\cap H_0^1(G)$
and $U=L^2(G_0)$.

  \item Let $\a(\cd)$ be a random variable with
Laplace distribution. Then, the system
\eqref{sch-system0} is exactly controllable in
average with $H=L^2(G)$, $V=H^4(G)\cap H_0^1(G)$
and $U=L^2(G_0)$.

  \item Let $\a(\cd)$ be a random variable with
standard Chi-squared distribution. Then, the
system \eqref{sch-system0} is exactly
controllable in average with $H=L^2(G)$,
$V=H^k(G)\cap H_0^1(G)$ and $U=L^2(G_0)$.
\end{itemize}
\end{theorem}

{\it Proof of the first conclusion in Theorem
\ref{con-th6}}\,: We divide the proof into two
steps.

{\bf Step 1}. In this step, we show that the
average of solutions can be represented by the
difference of the solutions of two Schr\"odinger
equations.  This was already noticed in
\cite{Z5}. We present  the argument here for the
sake of completeness.

We let $z_0\in L^2(G)$ in \eqref{sch-system1}
and $\tilde z(x,t;z_0)=\int_\O
z(x,T-t,\o;z_0)d\dbP(\o)$. Assume that
$z_0=\sum_{j=1}^\infty z_{0,j}e_j$. Similar to
\eqref{avu1}, we have that
\begin{equation}\label{savu1}
\tilde z(x,t;z_0) =
\frac{1}{b-a}\sum_{j=1}^\infty \frac{1}{i\l_j
t}z_{0,j}\big(e^{-i\l_j a t} - e^{-i\l_j b
t}\big)e_j.
\end{equation}
Let
$$
z_a(\cd,\cd)=\sum_{j=1}^\infty
\frac{z_{0,j}}{i\l_j}e^{-i\l_j a t},\q
z_b(\cd,\cd)=\sum_{j=1}^\infty
\frac{z_{0,j}}{i\l_j}e^{-i\l_j b t}.
$$
Clearly,  $z_a(\cd,\cd)$ and $z_b(\cd,\cd)$
solve the following equations respectively:
\begin{equation}\label{6.14-eq3}
\left\{
\begin{array}{ll}\ds
z_{a,t} + ai\D z_a =0 &\mbox{ in } G\times
[0,T],\\
\ns\ds z_a = 0 &\mbox{ on } \pa G\times
[0,T],\\
\ns\ds z_a(0)=A_\D^{-1}z_0 &\mbox{ in }G,
\end{array}
\right.
\end{equation}
\begin{equation}\label{6.14-eq4}
\left\{
\begin{array}{ll}\ds
z_{b,t} + bi\D z_b =0 &\mbox{ in } G\times
[0,T],\\
\ns\ds z_b = 0 &\mbox{ on } \pa G\times
[0,T],\\
\ns\ds z_b(0)=A_\D^{-1}z_0 &\mbox{ in }G.
\end{array}
\right.
\end{equation}
Further,
$$
t \tilde z(\cd,\cd) =  (z_a - z_b) (\cd,\cd).
$$

{\bf Step 2}. In this step, we establish the
exactly averaged observability estimate. From
\eqref{6.14-eq3} and \eqref{6.14-eq4}, we have
that
$$
\begin{array}{ll}\ds
\q(i\pa_t + a\D)(i\pa_t + b\D)(z_a -
z_b)\\
\ns\ds = (i\pa_t + a\D)(i\pa_t + b\D)
z_a - (i\pa_t + a\D)(i\pa_t + b\D) z_b\\
\ns\ds = (i\pa_t + b\D)(i\pa_t + a\D) z_a=0.
\end{array}
$$
Hence, we know that $(i\pa_t + b\D)(z_a-z_b)$
solves
\begin{equation}\label{6.20-eq1}
\left\{
\begin{array}{ll}\ds
i\f_t + a\D\f = 0 &\mbox{ in } G\times
(0,T],\\
\ns\ds \f = 0 &\mbox{ on } \pa G\times
(0,T),\\
\ns\ds \f(0) = (b-a) z_0 &\mbox{ in } G.
\end{array}
\right.
\end{equation}
By assumption ({\bf A3}),  for any $z_0\in
L^2(G)$, it holds that
\begin{equation}\label{6.14-eq5}
\begin{array}{ll}\ds
|z_0|_{L^2(G)}^2 \3n&\ds\leq C\int_0^T\int_{G_0}
|\f(x,t)|^2dxdt \leq C\int_0^T\int_{G_0}
|(i\pa_t +
b\D)(z_a-z_b)(x,t)|^2dxdt\\
\ns&\ds \leq C\int_0^T\int_{G_0} \big|(i\pa_t +
b\D)\big[t\tilde z(x,t;z_0)\big]\big|^2dxdt \leq
C\int_0^T\big|t\tilde
z(x,t;z_0)\big|^2_{H^2(G_0)}dt\\
\ns&\ds \leq C\int_0^T|\tilde
z(x,t;z_0)|^2_{H^2(G_0)}dt.
\end{array}
\end{equation}
\endpf

The proofs of the second to the fourth
conclusion in Theorem \ref{con-th6} are very
similar. We only give that for the second one.

\vspace{0.1cm}

{\it Proof of the second conclusion in Theorem
\ref{con-th6}}\,: Let $z_0\in V'=[H^2(\O)\cap
H_0^1(\O)]'$ in \eqref{sch-system1} and $\tilde
z(x,t;z_0)=\int_\O z(x,T-t,\o;z_0)d\dbP(\o)$. We
only need to prove that \vspace{-0.2cm}
\begin{equation}\label{6.16-eq1}
|z_0|_{V'}^2 \leq C\int_0^T\int_{G_0} |\tilde
z(x,t;z_0)|^2dxdt.
\end{equation}
We divide the proof into two steps.

{\bf Step 1}. In this step, we prove that a
``weak" version of the exact averaged
observability, that is, there is a lower order
term in the right hand side of the inequality.

Assume that $z_0=\sum_{j=1}^\infty z_{0,j}e_j\in
V'$. Then,\vspace{-0.2cm}
$$
z(x,t,\o;z_0) =  \sum_{j=1}^\infty
z_{0,j}e^{-i\a \l_j (T-t)}e_j.
$$
Similar to \eqref{ave1}, we have that
\vspace{-0.2cm}
\begin{equation}\label{save1}
\tilde z(x,t;z_0) = \int_\O
z(x,T-t,\o;z_0)d\dbP(\o) = \sum_{j=1}^\infty
\frac{1}{i\l_j t+1}z_{0,j} e^{-i\l_jt}e_j.
\end{equation}
This implies that for any $\d>0$,
\begin{equation}\label{10.17-eq1}
|\tilde z(\cd,\cd;z_0)|_{L^2(\d,T;L^2(G))}\leq
C(\d,T)|z_0|_{V'}.
\end{equation}

Let
$$
v(x,t)=\sum_{j=1}^\infty \frac{1}{i\l_j
t}z_{0,j} e^{-i\l_jt}e_j.
$$
From assumption ({\bf A3}), for a fixed $\d>0$,
we have that
\begin{equation}\label{5.31-eq1}
\begin{array}{ll}\ds
|z_0|_{V'}^2 \3n&\ds= \Big|\sum_{j=1}^\infty
\frac{1}{i\l_j}z_{0,j}
e^{-i\l_jt}e_j\Big|_{L^2(G)}^2 \leq
C\int_\d^T\int_{G_0}\Big|\sum_{j=1}^\infty
\frac{z_{0,j}}{\l_j}e^{-i\l_jt}e_j\Big|^2dxdt\\
\ns&\ds = C\int_\d^T\int_{G_0} |tv(x,t)|^2dxdt.
\end{array}
\end{equation}
Therefore,
\begin{equation}\label{5.31-eq2}
\begin{array}{ll}\ds
\q|z_0|_{V'}^2\\
\ns \ds \leq
C\int_\d^T\int_{G_0} |tv(x,t)|^2dxdt \\
\ns \ds\leq C\[\int_\d^T\int_{G_0} |t\tilde
z(x,t;z_0)|^2dxdt + \int_\d^T\int_{G_0}
|tv(x,t)-t\tilde z(x,t;z_0)|^2dxdt \]\\
\ns\ds \leq C \int_\d^T\int_{G_0} |t\tilde
z(x,t;z_0)|^2dxdt + C\int_\d^T\int_{G_0}
\Big|\sum_{j=1}^\infty \frac{1}{i\l_j(i\l_j
t+1)}z_{0,j} e^{-i\l_jt}e_j\Big|^2dxdt\\
\ns\ds \leq C \int_\d^T\int_{G_0} |t\tilde
z(x,t;z_0)|^2dxdt + C\int_\d^T\int_{G}
\Big|\sum_{j=1}^\infty \frac{1}{i\l_j(i\l_j
t+1)}z_{0,j} e^{-i\l_jt}e_j\Big|^2dxdt\\
\ns\ds \leq C \int_\d^T\int_{G_0} |t\tilde
z(x,t;z_0)|^2dxdt + C\int_\d^T\int_{G}
\Big|\sum_{j=1}^\infty \frac{1}{\l_j^2}z_{0,j}
e^{-i\l_jt}e_j\Big|^2dxdt\\
\ns\ds \leq C \int_\d^T\int_{G_0} |t\tilde
z(x,t;z_0)|^2dxdt + C |A_\D^{-1}z_0|_{V'}^2.
\end{array}
\end{equation}

{\bf Step 2}. In this step, we get rid of the
term $|A_\D^{-1}z_0|_{V'}^2$ in the right hand
side of \eqref{5.31-eq2} by a
compactness--uniqueness argument. More
precisely, we are going to prove that
\begin{equation}\label{10.20-eq7}
|z_0|_{V'}^2 \leq C \int_\d^T\int_{G_0} |t\tilde
z(x,t;z_0)|^2dxdt.
\end{equation}
If \eqref{10.20-eq7} is not true, then we can
find a sequence $\{z_0^n\}_{n=1}^\infty\subset
L^2(G)$ with $|z_0^n|_{H^{-2}(G)}=1$ such that
\begin{equation}\label{6.4-eq1}
\int_{\d}^T\int_{G_0} |t\tilde
z(x,t;z_0^n)|^2dxdt \leq \frac{1}{n}.
\end{equation}
Since $\{z_0^n\}_{n=1}^\infty$ is bounded in
$V'$, we can find a subsequence
$\{z_0^{n_k}\}_{k=1}^\infty\subset
\{z_0^n\}_{n=1}^\infty$ such that
\begin{equation}\label{10.24-eq1}
z_0^{n_k}  \mbox{ converges weakly to some
}z_0^*\in H^{-2}(G)
\end{equation}
and
\begin{equation}\label{10.17-eq3}
A_\D^{-1}z_0^{n_k} \mbox{ converges strongly to
} A_\D^{-1}z_0^* \mbox{ in }V'.
\end{equation}
According to \eqref{6.4-eq1} and
\eqref{5.31-eq2}, we know that
$$
|A_\D^{-1}z_0^{n_k}|_{V'}^2\geq \frac{1}{C} -
\frac{1}{n_k}.
$$
This, together with \eqref{10.17-eq3}, implies
that there is a positive constant $C>0$ such
that
\begin{equation}\label{10.24-eq2}
|A_\D^{-1}z_0^{*}|_{V'}^2\geq \frac{1}{C}.
\end{equation}
Thus, the limit $z_0^*$ is non trivial.

Further, from \eqref{save1} and
\eqref{10.17-eq1}, we know that $\tilde
z(\cd,\cd;z_0^{n_k})$ converges weakly to
$\tilde z(\cd,\cd;z_0^*)$ in $L^2(\d,T;L^2(G))$.
Hence,
$$
\int_{\d}^T\int_{G_0} |\tilde
z(x,t;z_0^*)|^2dxdt \leq
\lim_{k\to\infty}\int_{\d}^T\int_{G_0} |\tilde
z(x,t;z_0^{n_k})|^2dxdt \leq
\lim_{k\to\infty}\frac{1}{\d n_k}=0.
$$
Therefore, we find that
\begin{equation}\label{6.4-eq3}
\tilde z(\cd,\cd;z_0^*) =0 \;\mbox{ in }
G_0\times (\d,T).
\end{equation}

We would like to show that this leads to $z_0^*
\equiv 0$ which would then yield to a
contradiction.
To do this, we introduce the linear subspace
$$
\begin{array}{ll}\ds
\cE\=\big\{z_0\in V':\,\mbox{ the solution to
\eqref{sch-system1} with the initial datum $z_0$
fulfills } \\
\ns\ds \hspace{2.9cm} \tilde z(\cd,\cd;z_0) =0
\;\mbox{ in } G_0\times (\d,T) \big\}.
\end{array}
$$
Clearly, $z_0^*$ given in \eqref{10.24-eq1}
belongs to $\cE$. We want to prove that
$\cE=\{0\}$, which would be in contradiction
with the fact that $z_0^*$ is nonzero.

{\bf Step 3.} To show the claim that $\cE=\{0\}$
and conclude the proof, we proceed in several
steps.

{\bf Step 3.1.} We first prove that $\cE\subset
L^2(G)$.

To do this, given $\e\in (0,\d)$ and any
solution $\tilde z$, we introduce the discrete
time-derivative
\begin{equation}\label{6.20-eq3}
\hat z_\e(x,t;z_0)=\frac{\tilde
z(x,t+\e;z_0)-\tilde z(x,t;z_0)}{\e} \q\mbox{
for }\;t\in [0,T-\d].
\end{equation}
Then we have that
\begin{equation}\label{6.20-eq2}
\begin{array}{ll}\ds
\q\hat z_\e(x,t;z_0)\\
\ns\ds = \frac{1}{\e}\sum_{j=1}^\infty
\frac{1}{i\l_j (t+\e)+1}z_{0,j}
e^{-i\l_j(t+\e)}e_j -
\frac{1}{\e}\sum_{j=1}^\infty \frac{1}{i\l_j
t+1}z_{0,j} e^{-i\l_jt}e_j\\
\ns \ds = \frac{1}{\e}\!\sum_{j=1}^\infty
\frac{1}{i\l_j (t\!+\!\e)\!+\!1}z_{0,j}
\big(e^{-i\l_j(t+\e)}\! -\!
e^{-i\l_jt}\big)e_j\! +
\!\frac{1}{\e}\!\sum_{j=1}^\infty
\[\frac{1}{i\l_j (t\!+\!\e)\!+\!1}\!-\!\frac{1}{i\l_j t\!+\!1}\]z_{0,j}
e^{-i\l_jt}e_j \\
\ns \ds =  \sum_{j=1}^\infty \frac{1}{i\l_j
(t\!+\!\e)\!+\!\!1}z_{0,j}\frac{e^{-i\l_j\e}\!\!-\!1}{\e}
e^{-i\l_jt} e_j- \sum_{j=1}^\infty
\frac{i\l_j}{[i\l_j (t\!+\!\e)\!+\!1](i\l_j
t+1)}z_{0,j} e^{-i\l_jt}e_j
\end{array}
\end{equation}
and
\begin{equation}\label{6.20-eq2.1}
\hat z_\e(x,0;z_0)  =  \sum_{j=1}^\infty
\frac{1}{i\l_j
\e+1}z_{0,j}\frac{e^{-i\l_j\e}-1}{\e} e_j-
\sum_{j=1}^\infty \frac{i\l_j}{i\l_j
\e+1}z_{0,j} e_j.
\end{equation}
Let
\begin{equation}\label{6.20-eq4}
v_\e(x,t)= \sum_{j=1}^\infty \frac{1}{i\l_j
(t+\e)}z_{0,j}\frac{e^{-i\l_j\e}-1}{\e}
e^{-i\l_jt} e_j-\sum_{j=1}^\infty \frac{1}{i\l_j
(t\!+\!\e)}z_{0,j} e^{-i\l_jt}e_j.
\end{equation}
Then, again,  by assumption ({\bf A3}), as in
the proof of \eqref{5.31-eq1}, we have

$$
\Big|\frac{\tilde z_\e(x,\e;z_0)\!-\!\tilde
z_\e(x,0;z_0)}{\e}\Big|_{V'}^2\!\! = \big|\hat
z_\e(x,0;z_0)\big|_{V'}^2 \!\leq
C\int_\d^{T-\d}\!\!\!
\int_{G_0}\!|(t+\e)v_\e(x,t)|^2 dxdt,
$$
where the constant $C$ is independent of $\e$.
Thus, we obtain that
\begin{equation}\label{6.20-eq5}
\begin{array}{ll}\ds
\q\Big|\frac{\tilde z_\e(x,\e;z_0)-\tilde
z_\e(x,0;z_0)}{\e}\Big|_{V'}^2\\
\ns\ds\leq C\int_\d^{T-\d} \int_{G_0}|\hat
z_\e(x,t;z_0)|^2 dxdt + C\int_\d^{T-\d}
\int_{G_0}|\hat z_\e(x,t;z_0)-v_\e(x,t)|^2 dxdt
\\
\ns \ds \leq C\int_\d^{T-\d} \int_{G_0}|\hat
z_\e(x,t;z_0)|^2 dxdt + C\int_\d^{T-\d}
\int_{G}|\hat z_\e(x,t;z_0)-v_\e(x,t)|^2 dxdt.
\end{array}
\end{equation}
Let us estimate the second term in the right
hand side of \eqref{6.20-eq5}. From
\eqref{6.20-eq2} and \eqref{6.20-eq4}, we have
that
\begin{equation}\label{11.4-eq1}
\begin{array}{ll}\ds
\q\int_\d^{T-\d} \int_{G}|\hat z_\e(x,t;z_0)-v_\e(x,t)|^2 dxdt \\
\ns\ds \leq 2\int_\d^{T-\d}
\int_{G}\Big|\sum_{j=1}^\infty \[\frac{1}{i\l_j
(t+\e)+1}-\frac{1}{i\l_j
(t+\e)}\]z_{0,j}\frac{e^{-i\l_j\e}-1}{\e}
e^{-i\l_jt} e_j \Big|^2dxdt\\
\ns\ds \q + 2\int_\d^{T-\d}
\int_{G}\Big|\sum_{j=1}^\infty
\[\frac{i\l_j}{[i\l_j (t+\e)+1](i\l_j
t+1)}-\frac{1}{i\l_j (t+\e)}\]z_{0,j}
e^{-i\l_jt}e_j \Big|^2dxdt\\
\ns\ds \leq
C\sum_{j=1}^\infty\frac{1}{\l_j^4}z^2_{0,j}\Big|\frac{e^{-i\l_j\e}-1}{\e}\Big|^2
+ C\sum_{j=1}^\infty\frac{1}{\l_j^2}z^2_{0,j}.
\end{array}
\end{equation}
If $\l_j<\frac{1}{\e}$, then
\begin{equation}\label{11.4-eq2}
\Big|\frac{e^{-i\l_j\e}-1}{\e}\Big|=\frac{1}{\e}\Big|\sum_{k=1}^\infty
\frac{1}{k!}(-i\l_j\e)^k\Big|=\Big|i\l_j\sum_{k=1}^\infty
\frac{1}{k!}(-i\l_j\e)^{k-1} \Big|\leq
\l_j\sum_{k=1}^\infty\frac{1}{k!}\leq 2\l_j.
\end{equation}
If $\l_j>\frac{1}{\e}$, then
\begin{equation}\label{11.4-eq3}
\Big|\frac{e^{-i\l_j\e}-1}{\e}\Big|=\frac{1}{\e}\big|
e^{-i\l_j\e}-1\big|\leq \l_j\big|
e^{-i\l_j\e}-1\big|\leq 2\l_j.
\end{equation}
According to \eqref{11.4-eq1}, \eqref{11.4-eq2}
and \eqref{11.4-eq3}, we obtain that
\begin{equation}\label{11.4-eq4}
\int_\d^{T-\d} \int_{G}|\hat z_\e(x,t;z_0)-v_\e(x,t)|^2 dxdt \\
 \leq
C\sum_{j=1}^\infty\frac{1}{\l_j^2}z^2_{0,j}.
\end{equation}
As a result of \eqref{6.20-eq5} and
\eqref{11.4-eq4}, we have that
\begin{equation}\label{11.4-eq5}
\Big|\frac{\tilde z_\e(x,\e;z_0)-\tilde
z_\e(x,0;z_0)}{\e}\Big|_{V'}^2\leq
C\int_\d^{T-\d} \int_{G_0}|\hat z_\e(x,t;z_0)|^2
dxdt +
C\sum_{j=1}^\infty\frac{1}{\l_j^2}z^2_{0,j}.
\end{equation}
Since $z_0\in \cE$,  we know that
$$
\hat z_\e(x,t;z_0) = 0 \q\mbox{ in } G_0\times
(\d,T).
$$
This, together with \eqref{11.4-eq5}, implies
that for any $\e\in (0,\d)$,
\begin{equation}\label{6.20-eq6}
\Big|\frac{\tilde z_\e(x,\e;z_0)-\tilde
z_\e(x,0;z_0)}{\e}\Big|_{V'}^2\leq
C\sum_{j=1}^\infty \frac{z_{0,j}^2}{\l_j^2}\leq
C|z_0|_{V'}^2.
\end{equation}
Letting $\e\to 0$ in \eqref{6.20-eq6}, we obtain
that
\begin{equation}\label{6.20-eq6.1}
|z_0|_{L^2(G)}^2\leq C\Big|\sum_{j=1}^\infty
i\l_jz_{0,j} e_j\Big|_{V'}^2\leq
C\sum_{j=1}^\infty \frac{z_{0,j}^2}{\l_j^2}\leq
C|z_0|_{V'}^2.
\end{equation}
Then, we have that $\cE\subset L^2(G)$.

{\bf Step 3.2.} The same estimate deduces that
$\cE$ is a finite-dimensional subspace of
$L^2(G)$, since
\begin{equation}
|z_0|_{L^2(G)}^2\leq C|z_0|_{V'}^2.
\end{equation}

{\bf Step 3.3.} We now claim that
$A_\D\cE\subset \cE$.

Utilizing \eqref{save1} again, we have that
$$
\tilde z(\cd,\cd;z_0) \in C([0,T];L^2(G)) \mbox{
for any } z_0\in\cE.
$$
Since $\tilde z(\cd,\cd;z_0) =0$  in $G_0\times
(\d,T)$ for any $\d>0$, we see that $z_0=0$ in
$G_0$.

Thanks to $\cE\subset L^2(G)$, we get that
$A_\D\cE\subset V'$ and $\tilde z(x,t;A_\D
z_0)=A_\D\tilde z(x,t;z_0)=0$ in $G_0\times
(\d,T)$. Therefore, we obtain that
$A_\D\cE\subset \cE$.

{\bf Step 3.4.} To conclude, assume that
$\cE\neq \{0\}$. Then, there would exist a
non-trivial eigenfunction $\psi\in \cE$ and an
eigenvalue $\mu\in\dbR$  such that
$$
\left\{
\begin{array}{ll}\ds
-\D\psi=\mu\psi &\mbox{ in } G,\\
\ns\ds \psi = 0 &\mbox{ on } \pa G,\\
\ns\ds \psi = 0 &\mbox{ in } G_0.
\end{array}
\right.
$$
However, by the classical unique continuation
property for elliptic equations, $\psi=0$ in
$G$, which contradicts that $\psi$ is an
eigenfunction.

This concludes the proof of the fact  that
$\cE=\{0\}$.

As a consequence, we derive the desired estimate
\eqref{6.16-eq1}.
\endpf

\begin{remark}
We have utilized a compactness-uniqueness
argument in the above proof, which has been
extensively used in the proof of observability
estimates (see \cite{BLR} for example). Note
however that, normally, this is done for
solutions of PDE models. The averages under
consideration not being solutions of a specific
PDE this argument needs to be carefully adapted
as we have done above.
\end{remark}
%


\section{Further comments and open
problems}\label{sec-com}


\subsection{Further comments}


\subsubsection{Boundary averaged control for random heat
equations}


\q\;\,In this subsection, for convenience, we
assume that $\pa G$ is $C^\infty$ smooth,
although most comments and results make sense
with different weaker regularity assumptions.

We have solved the internal averaged
controllability problems for some particular
classes of random heat and Schr\"odinger
equations. The same could be done for boundary
control problems.

Let us consider the following  heat equation
with boundary control and random constant
diffusivity:
\begin{equation}\label{heat-system4}
\left\{
\begin{array}{ll}\ds
y_t - \a \D y =0 &\mbox{ in }G\times
(0,T],\\
\ns\ds y=u &\mbox{ on } \G_0\times (0,T),\\
\ns\ds y=0 &\mbox{ on } (\pa
G\setminus\G_0)\times (0,T),
\\
\ns\ds y(0)=y_0 &\mbox{ in }G.
\end{array}
\right.
\end{equation}
Here $\G_0$ is an open subset of $\pa G$,
$\a(\cd)$ is a random variable, $u\in
L^2(0,T;L^2(\G_0))$ and $y_0\in L^2(G)$.

We have the following result.

\begin{theorem}\label{con-th7}
Let $\a(\cd)$ be a uniformly distributed or an
exponentially distributed random variable. The
system \eqref{heat-system4} is null controllable
in average with control $u(\cd)\in
L^2(0,T;L^2(\G_0))$. Further, there is a
constant $C>0$ such that
\begin{equation}\label{th7-eq1}
|u|_{L^2(0,T;L^2(\G_0))} \leq C|y_0|_{L^2(G)}.
\end{equation}
\end{theorem}

To prove Theorem \ref{con-th7}, we consider the
adjoint system of \eqref{heat-system4} as
follows:
\begin{equation}\label{heat-system5}
\left\{
\begin{array}{ll}\ds
z_t + \a \D z =0 &\mbox{ in }G\times
(0,T],\\
\ns\ds z=0 &\mbox{ on } \pa G\times (0,T),\\
\ns\ds z(T)=z_0 &\mbox{ in }G.
\end{array}
\right.
\end{equation}
One only need to prove the following result.
\begin{theorem}\label{th4}
There exists a constant $C>0$ such that for any
$z_0\in L^2(G)$, and either $\a(\cd)$ is a
uniformly distributed or an exponentially
distributed random variable, it holds that
\begin{equation}\label{th4-eq1}
\Big|\int_\O
z(\cd,0,\o;z_0)d\dbP(\o)\Big|_{L^2(G)}^2 \leq
C\int_0^T\int_{\G_0} \Big|\int_\O \frac{\pa
z(x,t,\o;z_0)}{\pa\nu}d\dbP(\o)\Big|^2d\G_0dt.
\end{equation}
\end{theorem}
The proof is very similar to the one for Theorem
\ref{th0}. We only give a sketch here. Let us
assume that $\a(\cd)$ is an exponentially
distributed random variable. From \cite[Page
345]{LR}, we have the following result:
\begin{equation}\label{6.21-eq1}
\sum_{\l_j\leq r} a_j^2 \leq
C_2e^{C_2\sqrt{r}}e^{\frac{1}{t_2-t_1}}\int_{t_1}^{t_2}\int_{\G_0}\Big|\sum_{\l_j\leq
r}e^{\sqrt{\l_j}t}a_j \frac{\pa
e_j}{\pa\nu}\Big|^2 d\G_0dt \q\mbox{ for any }
0\leq t_1<t_2\leq T.
\end{equation}

Let $T_k = (1-\frac{1}{2^{k-1}})T$ for
$k\in\dbN$ and $r_k=2^{2(k+1)}[\ln(6C_2)+C_2]$.
Similar to the proof of \eqref{5.28-eq10}, we
can obtain that
\begin{equation}\label{6.21-eq2}
\begin{array}{ll}\ds
\frac{T_{k+1}-T_k}{6C_2e^{C_2\sqrt{r_k}}}|\tilde
z(\cd,T_{k};z_0)|_{L^2(G)} -
\frac{C_2e^{C_2\sqrt{r_k}}+1}{C_2e^{C_2\sqrt{r_k}}}
(T_{k+1}-T_k)e^{-
\frac{T_{k+1}-T_k}{2}r_k}|\tilde
z(\cd,T_{k+1};z_0)|_{L^2(G)}\\
\ns\ds \leq \int_{T_k}^{T_{k+1}}\int_{\G_0}
\Big|\frac{\pa\tilde
z(x,t;z_0)}{\pa\nu}\Big|^2d\G_0 dt.
\end{array}
\end{equation}
By summarizing the inequality \eqref{6.21-eq2}
from $k=1$ to $k=\infty$, we obtain that
\begin{equation}\label{5.28-eq11.1}
\frac{T_{2}-T_{1}}{6C_2e^{C_2\sqrt{r_{1}}}}
|z_0|_{L^2(G)} + \sum_{k=1}^\infty f_k |\tilde
z(\cd,T_{k+1};z_0)|_{L^2(G)}\leq \int_0^T
\int_{\G_0}\Big|\frac{\pa\tilde
z(x,t;z_0)}{\pa\nu}\Big|d\G_0dt,
\end{equation}
where
$$
f_k=\frac{T_{k+2}-T_{k+1}}{6C_2e^{C_2\sqrt{r_{k+1}}}}-
\frac{C_2e^{C_2\sqrt{r_k}}+1}{C_2e^{C_2\sqrt{r_k}}}(T_{k+1}-T_k)e^{-
\frac{T_{k+1}-T_k}{2}r_k}, \q k=1,2,\cds
$$
From  \eqref{5.28-eq3} and $r_k
=2^{2(k+1)}[\ln(6C_2)+C_2]$, we have that
$$
f_k\geq 0 \;\mbox{ for any }\; k=1, 2, \cds
$$
This, together with \eqref{5.28-eq11.1}, deduces
that
\begin{equation}\label{6.25-eq2}
|z_0|_{L^2(G)}^2\leq
\frac{6C_2e^{C_2\sqrt{r_{1}}}}{T_{2}-T_{1}}\int_0^T
\int_{\G_0}\Big|\frac{\pa \tilde
z(x,t;z_0)}{\pa\nu}\Big|^2d\G_0dt.
\end{equation}
This completes the proof.
\endpf


\subsubsection{Boundary averaged control
for random Schr\"odinger equations}


Consider the following random  Schr\"odinger
equations:
\begin{equation}\label{sch-system3}
\left\{
\begin{array}{ll}\ds
y_t - \a i\D y =0 &\mbox{ in }G\times
(0,T],\\
\ns\ds y=u &\mbox{ on } \G_0\times (0,T),\\
\ns\ds y=0 &\mbox{ on } (\pa
G\setminus\G_0)\times (0,T),
\\
\ns\ds y(0)=y_0 &\mbox{ in }G.
\end{array}
\right.
\end{equation}
Here $\G_0$ is an open subset of $\pa G$,
$\a(\cd)$ is a random variable, $u\in
L^2(0,T;L^2(\G_0))$ and $y_0\in L^2(G)$. We have
the following controllability results.

\begin{theorem}\label{con-th8}
System \eqref{sch-system3} is null controllable
in average if $\a$ is a random variable with
normal distribution or Cauchy distribution.
\end{theorem}

Further, we assume the following condition
holds:

({\bf A4}) Whatever $T>0$ and $k \ne 0$ are,
there are constants $C_3$ and $C_4$ such that
for any $\f_0\in H_0^1(G)$, the solution
$\f(\cd,\cd)$ to \eqref{Sch-system2} satisfies
\begin{equation}\label{6.29-eq2.1}
|\f_0|_{H_0^1(G)}^2\leq C_3\int_0^T\int_{G_0}
\Big|\frac{\pa\f}{\pa\nu}\Big|^2d\G_0 dt\leq
C_4|\f_0|_{H_0^1(G)}^2.
\end{equation}
We refer the readers to \cite{Le} for the
conditions on $\G_0$ for which ({\bf A4}) holds.

\begin{theorem}\label{con-th6.1}
Assume that ({\bf A4}) holds. We have the
following results:
\begin{itemize}
  \item  Let $\a(\cd)$
be a uniformly distributed random variable on an
interval $[a,b]$. Then, the system
\eqref{sch-system0} is exactly controllable in
average with $V=H=L^2(G)$ and $U=H^{-1}(\G_0)$.

  \item  Let $\a(\cd)$ be an exponentially
distributed random variable. Then, the system
\eqref{sch-system0} is exactly  controllable in
average with $H=L^2(G)$, $V=H^1_0(G)$ and
$U=L^2(\G_0)$.

\item  Let $\a(\cd)$ be a random variable with
Laplace distribution. Then, the system
\eqref{sch-system0} is exactly  controllable in
average with $H=L^2(G)$, $V=H^3(G)\cap H_0^1(G)$
and $U=L^2(\G_0)$.

\item  Let $\a(\cd)$ be a random variable with
standard Chi-squared distribution. Then, the
system \eqref{sch-system0} is exactly
controllable in average with $H=L^2(G)$,
$V=H^{k-1}(G)\cap H_0^1(G)$ and $U=L^2(\G_0)$.
\end{itemize}
\end{theorem}

As usual, we introduce the adjoint system of
\eqref{sch-system3} as follows:
\begin{equation}\label{sch-system4}
\left\{
\begin{array}{ll}\ds
z_t + i\a\D z = 0 &\mbox{ in } G\times
(0,T),\\
\ns\ds z=0 &\mbox{ on } \pa G\times
(0,T),\\
\ns\ds z(T)=z_0 &\mbox{ in }G.
\end{array}
\right.
\end{equation}

We can prove the following results.
\begin{theorem}\label{ob-th2}
Let  $\a(\cd)$ be a standard normally
distributed random variable. There exists a
constant $C>0$ such that for any $z_0\in
L^2(G)$, it holds that
\begin{equation}\label{ob-th2-eq1}
\Big|\sum_{j=1}^\infty z_{0,j}e^{-\l_j^2
T^2}e_j\Big|_{L^2(G)}^2 \leq
C\int_0^T\int_{\G_0} \Big|\int_\O\frac{\pa
z(x,t,\o;z_0)}{\pa\nu}d\dbP(\o)\Big|^2d\G_0dt.
\end{equation}
\end{theorem}
\begin{theorem}\label{ob-th3}
Let $\a(\cd)$ be a random variable with standard
Cauchy distribution. There exists a positive
constant $C$  such that for any $z_0\in L^2(G)$,
it holds that
\begin{equation}\label{ob-th3-eq1}
\Big|\sum_{j=1}^\infty z_{0,j}e^{-
\l_jT}e_j\Big|_{L^2(G)}^2 \leq C\int_0^T
\int_{\G_0} \Big|\int_\O\frac{\pa
z(x,t,\o;z_0)}{\pa\nu}d\dbP(\o)\Big|^2d\G_0dt.
\end{equation}
\end{theorem}

By means of \eqref{savc1}, Theorem \ref{ob-th3}
is nothing but the boundary observability
estimate for heat equations. The proof of
Theorem \ref{ob-th2} is very similar to the one
for Theorem \ref{th4}. We omit it here.

Further, the proof of Theorem \ref{con-th6.1} is
also very analogous to the proofs for Theorems
\ref{con-th6}. We only give a sketch of the
proof of the second conclusion in Theorem
\ref{con-th6.1}.

Let $z_0\in H^{-1}(G)$ in \eqref{sch-system1}
and $\tilde z(x,t;z_0)=\int_\O
z(x,T-t,\o;z_0)d\dbP(\o)$. Then, we only need to
prove that
\begin{equation}\label{6.16-eq1.1}
|z_0|_{H^{-1}(G)}^2 \leq C\int_0^T\int_{\G_0}
\Big|\frac{\pa\tilde
z(x,t;z_0)}{\pa\nu}\Big|^2d\G_0dt.
\end{equation}

Assume that $z_0=\sum_{j=1}^\infty z_{0,j}e_j$.
We have
\begin{equation}\label{6.21-eq3}
\tilde z(x,t;z_0)= \sum_{j=1}^\infty
\frac{1}{i\l_j t+1}z_{0,j} e^{-i\l_jt}e_j.
\end{equation}
Let
$$
v(x,t)=\sum_{j=1}^\infty \frac{1}{i\l_j
t}z_{0,j} e^{-i\l_jt}e_j.
$$
From ({\bf A4}), for a fixed $\d\in (0,T)$, we
have that
\begin{equation}\label{5.31-eq1.1}
|z_0|_{H^{-1}(G)}^2 = \Big|\sum_{j=1}^\infty
\frac{1}{i\l_j }z_{0,j}
e^{-i\l_jt}e_j\Big|_{H_0^1(G)}^2 \leq
C\int_\d^T\int_{\G_0} \Big|t\frac{\pa
v(x,t)}{\pa\nu}\Big|^2d\G_0dt.
\end{equation}
Therefore,
\begin{equation}\label{6.21-eq4}
\begin{array}{ll}\ds
\q|z_0|_{H^{-1}(G)}^2\\
\ns \ds \leq
C\int_\d^T\int_{\G_0} \Big|t\frac{\pa v(x,t)}{\pa\nu}\Big|^2 d\G_0dt \\
\ns \ds\leq C\[\int_\d^T\int_{\G_0}
\Big|t\frac{\pa\tilde
z(x,t;z_0)}{\pa\nu}\Big|^2d\G_0dt +
\int_\d^T\int_{\G_0}
\Big|t\frac{\pa v(x,t)}{\pa\nu}-t\frac{\pa\tilde z(x,t;z_0)}{\pa\nu}\Big|^2d\G_0dt \]\\
\ns\ds \leq C \int_\d^T\int_{\G_0}
\Big|t\frac{\pa\tilde
z(x,t;z_0)}{\pa\nu}\Big|^2d\G_0dt + C\int_\d^T
\Big|\sum_{j=1}^\infty \frac{1}{i\l_j(i\l_j
t+1)}z_{0,j} e^{-i\l_jt}e_j\Big|^2_{H^{1}(G)}dt \\
\ns\ds \leq C \int_\d^T\int_{\G_0}
\Big|t\frac{\pa\tilde
z(x,t;z_0)}{\pa\nu}\Big|^2d\G_0dt + C
|z_0|_{H^{-3}(G)}^2.
\end{array}
\end{equation}

We claim that
\begin{equation}\label{6.26-eq1}
|z_0|_{H^{-1}(G)}^2  \leq C \int_0^T\int_{\G_0}
\Big|t\frac{\pa\tilde
z(x,t;z_0)}{\pa\nu}\Big|^2d\G_0dt.
\end{equation}
If \eqref{6.26-eq1} is not true, then we can
find a sequence $\{z_0^n\}_{n=1}^\infty\subset
L^2(G)$ with $|z_0^n|_{H^{-2}(G)}=1$ such that
\begin{equation}\label{6.26-eq3}
\int_0^T\int_{\G_0} \Big|t\frac{\pa\tilde
z(x,t;z_0)}{\pa\nu}\Big|^2d\G_0dt \leq
\frac{1}{n}.
\end{equation}
Since $\{z_0^n\}_{n=1}^\infty$ is bounded in
$H^{-1}(G)$, we can find a subsequence
$\{z_0^{n_k}\}_{k=1}^\infty\subset
\{z_0^n\}_{n=1}^\infty$ such that $z_0^{n_k}$
converges weakly to some $z_0^*\in H^{-1}(G)$.
From \eqref{6.26-eq3}, we know that
\begin{equation}\label{6.21-eq6.1}
|z_0^*|_{H^{-3}(G)}^2\geq \frac{1}{C}
\end{equation}
for a positive constant and
\begin{equation}\label{6.21-eq6}
\frac{\pa\tilde z(\cd,\cd;z_0^*)}{\pa\nu} =0
\;\mbox{ on } \G_0\times (0,T).
\end{equation}
Put
$$
\wt\cE\=\{z_0\in H^{-1}(G):\,\mbox{ the solution
to \eqref{sch-system4} with the final datum
$z_0$ fulfills \eqref{6.21-eq6}} \}.
$$
Analogous to the proof that $\cE=\{0\}$, we can
prove that $\wt\cE=\{0\}$, which contradicts
\eqref{6.21-eq6.1}. Hence, we know that
\eqref{6.26-eq1} holds.
\endpf


\subsubsection{Averaged control for random heat
and Schr\"odinger equations from measurable
sets}


\q\,We have considered the averaged control
problems of random heat and Schr\"odinger
equations for the internal control (\resp
boundary control) supported in $G_0\times E$
(\resp $\G_0\times E$), where $G_0\subset G$
(\resp $\G_0\subset \pa G$) is a nonempty open
subset. By means of the method in \cite{Ap}, one
can consider the case that the internal control
(\resp the boundary control) is supported in a
measurable subset $\cD\subset G\times (0,T)$
(\resp $\cG\subset \pa G\times (0,T)$). For
instance, we can consider the following systems:
\begin{equation}\label{heat-system7}
\left\{
\begin{array}{ll}\ds
y_t - \a\D y = \chi_{\cD}u &\mbox{ in } G\times
(0,T),\\
\ns\ds y=0 &\mbox{ on } \pa G\times
(0,T),\\
\ns\ds y(0)=y_0 &\mbox{ in }G,
\end{array}
\right.
\end{equation}
and
\begin{equation}\label{sch-system6}
\left\{
\begin{array}{ll}\ds
y_t - i\a\D y = \chi_{\cD}u &\mbox{ in } G\times
(0,T),\\
\ns\ds y=0 &\mbox{ on } \pa G\times
(0,T),\\
\ns\ds y(0)=y_0 &\mbox{ in }G.
\end{array}
\right.
\end{equation}
Here $y_0 \in L^2(G)$, $\cD\subset G\times
(0,T)$ is a Lebesgue measurable set with
positive Lebesgue measure and $u\in
L^\infty(G\times (0,T))$. One can combine the
proof of Theorem \ref{con-th1} and Corollary 1
in \cite{Ap} to prove the following results.

\begin{theorem}\label{con-th9}
Assume that, either $\a(\cd)$ is a uniformly or
exponentially distributed random variable. Then,
the system \eqref{heat-system7} is null
controllable in average with control $u(\cd)\in
L^\infty(G\times (0,T))$. Further, there is a
constant $C>0$ such that
$$
|u|_{L^\infty(G\times (0,T))} \leq
C|y_0|_{L^2(G)}.
$$
\end{theorem}

\begin{theorem}\label{con-th10}
If $\a$ is a random variable with normal
distribution or Cauchy distribution, then the
system \eqref{sch-system6} is null  controllable
in average with control $u(\cd)\in
L^\infty(G\times (0,T))$. Further, there is a
constant $C>0$ such that
$$
|u|_{L^\infty(G\times (0,T))} \leq
C|y_0|_{L^2(G)}.
$$
\end{theorem}

Next, let us consider the following systems:
\begin{equation}\label{heat-system8}
\left\{
\begin{array}{ll}\ds
y_t - \a \D y =0 &\mbox{ in }G\times
(0,T],\\
\ns\ds y=u &\mbox{ on } \cG,\\
\ns\ds y=0 &\mbox{ on } [\pa G\times
(0,T)]\setminus\cG,
\\
\ns\ds y(0)=y_0 &\mbox{ in }G,
\end{array}
\right.
\end{equation}
and
\begin{equation}\label{sch-system7}
\left\{
\begin{array}{ll}\ds
y_t - \a i\D y =0 &\mbox{ in }G\times
(0,T],\\
\ns\ds y=u &\mbox{ on } \cG,\\
\ns\ds y=0 &\mbox{ on } [\pa G\times
(0,T)]\setminus\cG,
\\
\ns\ds y(0)=y_0 &\mbox{ in }G.
\end{array}
\right.
\end{equation}
Here $\cG$ is a Lebesgue measurable subset of
$\pa G\times(0,T)$ with positive Lebesgue
measure, $\a(\cd)$ is a random variable, $u\in
L^\infty(\pa G\times(0,T))$ and $y_0\in L^2(G)$.

One can combine the proof of Theorem
\ref{con-th7} and Corollary 1 in \cite{Ap} to
prove the following results.

\begin{theorem}\label{con-th11}
Assume that, either $\a(\cd)$ is a uniformly or
exponentially distributed random variable. Then,
the system \eqref{heat-system8} is null
controllable in average with control $u(\cd)\in
L^\infty(G\times (0,T))$. Further, there is a
constant $C>0$ such that
$$
|u|_{L^\infty(\pa G\times (0,T))} \leq
C|y_0|_{L^2(G)}.
$$
\end{theorem}

\begin{theorem}\label{con-th12}
If $\a$ is a random variable with normal
distribution or Cauchy distribution, then the
system \eqref{sch-system7} is null  controllable
in average with control $u(\cd)\in
L^\infty(G\times (0,T))$. Further, there is a
constant $C>0$ such that
$$
|u|_{L^\infty(\pa G\times (0,T))} \leq
C|y_0|_{L^2(G)}.
$$
\end{theorem}
%

\subsubsection{Internal averaged control
for random fractional Schr\"odinger equations}


\q\,We have studied the averaged controllability
problems for some random heat equations and
random Schr\"odinger equations. In the results
proved so far we have obtained averaged
controllability for parameter-depending
equations that were controllable for each value
of the parameter. Here, we give an example of
model which is not null controllable when one
fixes an $\o$ but that gains null
controllability by the averaging process.

Consider the following equation:
\begin{equation}\label{fsch-system0}
\left\{
\begin{array}{ll}\ds
iy_t + \a A_\D^{\g} y = B u &\mbox{ in } (0,T],\\
\ns\ds y(0)=y_0.
\end{array}
\right.
\end{equation}
Here $\g\in (\frac{1}{4},\frac{1}{2})$, $y_0\in
L^2(G)$, $u\in L^2(0,T;L^2(G_0))$ and
$Bu=\chi_{G_0}u$.

The adjoint system of \eqref{fsch-system0} reads
\begin{equation}\label{fsch-system1}
\left\{
\begin{array}{ll}\ds
iz_t - \a A_\D^{\g} z = 0 &\mbox{ in } [0,T),\\
\ns\ds z(T)=z_0,
\end{array}
\right.
\end{equation}
where $z_0\in L^2(G)$. Assume that
$z_0=\sum_{j=1}^\infty z_{0,j}e_j$. If
$\a(\cd):\O\to\dbR$ is a standard normally
distributed random variable, then we know that
\begin{equation}\label{6.22-eq1}
\begin{array}{ll}\ds
\int_\O z(x,t,\o;z_0)d\dbP(\o)\3n&\ds =
\frac{1}{\sqrt{2\pi}}\int_{-\infty}^\infty
e^{-\frac{\a^2}{2}}\sum_{j=1}^\infty
z_{0,j}e^{-i\l_j^\g \a
(T-t)}e_jd\a\\
\ns&\ds= \sum_{j=1}^\infty z_{0,j}e^{-\l_j^{2\g}
(T-t)^2}e_j.
\end{array}
\end{equation}
Similar to the proof of Theorem \ref{th2}, we
can establish the following result.
\begin{theorem}\label{ob-th4}
Let $E\subset [0,T]$ be a measurable set with
positive Lebesgue measure $m(E)$. There exists a
constant $C>0$ such that for any $z_0\in
L^2(G)$, it holds that
\begin{equation}\label{th2-eq1.1}
\Big|\sum_{j=1}^\infty z_{0,j}e^{-\l_j^{2\g}
T^2}e_j\Big|_{L^2(G)} \leq C\int_E \Big|\int_\O
z(x,t,\o;z_0)d\dbP(\o)\Big|_{L^2(G_0)}dt.
\end{equation}
\end{theorem}
Theorem \ref{ob-th4} implies that the system
\eqref{fsch-system0} is null  controllable in
average. However,  it is not null controllable
for any fixed $\a\in\dbR$ even for $d=1$. For
example, let $\a=1$ and $G=(0,1)$. Then the
solution to the system \eqref{fsch-system1}
reads
$$
z(x,t)=\sqrt{2}\sum_{j=1}^\infty
z_{0,j}e^{-i(j\pi)^{2\g}t}\sin(j\pi x).
$$
Since
$$
\lim_{j\to\infty}
\big\{[(j+1)\pi)]^{2\g}-(j\pi)^{2\g}\big\}=0,
$$
we know that the following inequality
$$
|z(\cd,T)|_{L^2(0,1)}^2 = \sum_{j=1}^\infty
z_{0,j}^2 \leq C\int_0^T\int_{G_0}|z(x,t)|^2dxdt
$$
does not hold for any $C>0$.


\subsubsection{Internal averaged control
for random heat equations with variable
coefficients}


\q\,We can also consider the approximate
averaged controllability problem of some more
general random heat equations. We make the
following assumptions on the coefficients
$a^{jk}: \;\cl{G}\times \O \to \dbR^{n\times
n}\;$ ($j,k=1,2,\cdots,n$):

\medskip

{\bf (H1)} {\it $a^{jk}(\cd,\o):\;\cl{G}\to
\dbR$ is analytic, $\dbP$-a.s., and
$a^{jk}=a^{kj}$.}

\medskip

{\bf (H2)} {\it For a.e. $\o\in\O$, there is a
constant $C(\o)>0$ such that for any multi-index
$\eta=(\eta_1,\cds,\eta_n)\in
(\dbN\cup\{0\})^n$,
$$
\Big|\frac{\pa^\eta a^{jk}(x,\o)}{\pa
x^\eta}\Big|\leq
\frac{C(\o)|\eta|!}{R^{|\eta|}},\q\mbox{ for any
} j,k=1,\cds,n,
$$
where $R$ is a positive constant larger than
$\max_{x\in \cl{G}}|x|$, and $C(\cd)$ satisfies
that
$$
\int_\O C(\o)d\dbP(\o)<\infty.
$$
}

\medskip

{\bf (H3)} {\it There exists a constant $s_0>0$
such that
\begin{equation}\label{h1}
\sum_{j,k=1}^n a^{jk}(\o,t,x)\xi^{j}\xi^{k} \geq
s_0|\xi|^{2},\q \forall\,(\o,x,\xi)\equiv
(\o,x,\xi^{1},\cdots,\xi^{n}) \in \Omega\times G
 \times \dbR^{n}.
\end{equation}}

Consider the following random heat equation:
\begin{equation}\label{heat-system2}
\left\{
\begin{array}{ll}\ds
y_t - \sum_{j,k=1}^n
\big(a^{jk}y_{x_j}\big)_{x_k} = \chi_{G_0}u
&\mbox{ in } G\times
(0,T),\\
\ns\ds y=0 &\mbox{ on } \pa G\times
(0,T),\\
\ns\ds y(0)=y_0 &\mbox{ in }G.
\end{array}
\right.
\end{equation}
Here the initial datum $y_0 \in L^2(G)$.

We have the following result.
\begin{theorem}\label{th-ap-con2}
Under the assumptions {\bf (H1)}--{\bf (H3)}
above system \eqref{heat-system2} is
approximately controllable in average in any
time $T>0$ and from any open non-empty subset
$G_0$ of $G$.
\end{theorem}

{\it Proof}\,: By Theorem \ref{th-ap-ob}, we
only need to prove that the adjoint system of
\eqref{heat-system2} satisfies the averaged
unique continuation property. Its adjoint system
reads
\begin{equation}\label{heat-system3}
\left\{
\begin{array}{ll}\ds
z_t + \sum_{j,k=1}^n
\big(a^{jk}z_{x_j}\big)_{x_k} = 0 &\mbox{ in }
G\times
(0,T),\\
\ns\ds z=0 &\mbox{ on } \pa G\times
(0,T),\\
\ns\ds z(T)=z_0 &\mbox{ in }G,
\end{array}
\right.
\end{equation}
where the final datum $z_0 \in L^2(G)$. From
{\bf (H1)} and {\bf (H3)}, we know that for any
$t\in [0,T)$ and a.e. $\o\in\O$, the solution
$z(\cd,t,\o;z_0)$ is analytic in $G$ (see
\cite{F,KN} for example). Further, for any ball
$B_r \subset G$ with radius $r$, there is a
constant $C>0$ such that for any multi-index
$\eta \in (\dbN\cup\{0\})^n$,
$$
\Big|\frac{\pa^\eta z(\cd,t,\o;z_0)}{\pa
x^\eta}\Big|\leq
CC(\o)\frac{|\eta|!}{r^{|\eta|}} \;\mbox{ in }
B_r.
$$
From {\bf (H2)}, we have that
$$
\Big|\frac{\pa^\eta \int_\O
z(\cd,t,\o)d\dbP(\o)}{\pa x^\eta}\Big|\leq
C\frac{|\eta|!}{r^{|\eta|}}.
$$
Then, we know that $\int_\O
z(\cd,t,\o)d\dbP(\o)$ is analytic in $B_r$.
Hence, it is analytic in $G$. Since $\int_\O
z(\cd,\cd,\o)d\dbP(\o) = 0$ in $G_0\times
(0,T)$, we get that it vanishes everywhere in
$G\times (0,T)$. Noting that it is continuous in
$L^2(G)$ with respect to $t$, we conclude that
$z_0=0$ in $G$, which implies that
\eqref{heat-system3} satisfies the averaged
unique continuation property.
\endpf


\subsubsection{Averaged controllability problems
for the random heat and the random Schr\"odinger
equations random initial data}


\q\, One can also consider the internal and
boundary averaged controllability problems for
random heat and random Schr\"odinger equations
with random initial data. Let us first consider
the following random heat equation:
\begin{equation}\label{heat-system9}
\left\{
\begin{array}{ll}\ds
y_t - \a\D y = \chi_{G_0\times E}u &\mbox{ in }
G\times
(0,T),\\
\ns\ds y=0 &\mbox{ on } \pa G\times
(0,T),\\
\ns\ds y(0,\o)=y_0(\o) &\mbox{ in }G,
\end{array}
\right.
\end{equation}
Here $y_0(\cd) \in L^2(\O;L^2(G))$, and $G_0$
and $E$ are subsets of $G$ and $[0, T]$
respectively, where the controls are being
applied.  The constant diffusivity $\a:\O\to
\dbR^+$ is assumed to be a random variable.

According to Remark \ref{rmk2}, we know that to
prove the averaged null controllability of
\eqref{heat-system9}, we only need to establish
the following observability estimate:
\begin{equation}\label{10.20-eq1}
\(\int_\O\Big|z(0,\o;z_0)\Big|^2_{L^2(G)}d\dbP(\o)\)^{\frac{1}{2}}
\leq C\int_E \Big|\int_\O
z(t,\o;z_0)d\dbP(\o)\Big|_{L^2(G_0)}dt,
\end{equation}
where $z$ solves \eqref{heat-system1} and $C$ is
independent of $z_0$.

If $\a(\cd)$ is the exponentially distributed
random variable, then
\begin{equation}\label{10.20-eq2}
\begin{array}{ll}\ds
\int_\O\Big|z(0,\o;z_0)\Big|^2_{L^2(G)}d\dbP(\o)\3n&\ds=
\int_1^\infty e^{-(\a-1)}\sum_{j=1}^\infty
z^2_{0,j}e^{-2\l_j \a
T}d\a\\
\ns&\ds = \sum_{j=1}^\infty \frac{1}{2\l_j
T+1}z_{0,j}^2 e^{-2\l_jT}.
\end{array}
\end{equation}
From \eqref{10.20-eq2}, similar to the proof of
the inequality \eqref{th0-eq1}, one can obtain
\eqref{10.20-eq1}. The same thing can be done if
$\a(\cd)$ is a uniformly distributed random
variable on $[a,b]$ for $0<a<b$.

\begin{proposition}\label{con-prop1}
Let ({\bf A1}) and ({\bf A2}) hold. Assume that,
either $\a(\cd)$ is a uniformly or exponentially
distributed random variable. Then, the system
\eqref{heat-system9} is null  controllable in
average with control $u(\cd)\in
L^2(0,T;L^2(G_0))$. Further, there is a constant
$C>0$ such that
\begin{equation}\label{prop1-eq1}
|u|_{L^2(0,T;L^2(G_0))} \leq
C|y_0|_{L^2(\O;L^2(G))}.
\end{equation}
\end{proposition}

Thanks to Remark \ref{rmk2}, we know that in
order to show that \eqref{heat-system9} is
 approximately controllable in average, one just needs to
prove that the solution to \eqref{heat-system1}
satisfies the averaged unique continuation
property, which is obtained in the proof of
Theorem \ref{th-ap-con1}. Hence, we have the
following result:
\begin{proposition}\label{con-prop2}
Let ({\bf A1}) and ({\bf A2}) hold.  System
\eqref{heat-system9} is approximately
controllable in average, provided that $\a(\cd)$
is a uniformly distributed or an exponentially
distributed random variable.
\end{proposition}

Next, we consider the averaged controllability
problem for the following random \linebreak
Schr\"odinger equation:
\begin{equation}\label{sch-system8}
\left\{
\begin{array}{ll}\ds
y_t - i\a\D y = \chi_{G_0\times E}u &\mbox{ in }
G\times
(0,T),\\
\ns\ds y=0 &\mbox{ on } \pa G\times
(0,T),\\
\ns\ds y(0)=y_0 &\mbox{ in }G.
\end{array}
\right.
\end{equation}
Here $y_0(\cd) \in L^2(\O;V)$, and $G_0$ and $E$
are subsets of $G$ and $[0, T]$ respectively,
where the controls are being applied.  $\a:\O\to
\dbR$ is assumed to be a random variable.

In virtue of Remark \ref{rmk2}, we know that to
get the averaged exact controllability of the
system \eqref{sch-system8}, one only need to
prove that the solution to the equation
\eqref{sch-system1} is exactly averaged
observable. As a result of these facts, we know
that the conclusions in Theorem \ref{con-th6}
also hold for the system \eqref{sch-system8}.
More precisely, we have the following results:
\begin{proposition}\label{con-prop4}
The following results hold:
\begin{itemize}
  \item Let $\a(\cd)$ be a uniformly distributed
random variable on an interval $[a,b]$. Then,
the system \eqref{sch-system8} is exactly
controllable in average with $V=H=L^2(G)$ and
$U=H^{-2}(G_0)$.

  \item Let $\a(\cd)$ be an exponentially
distributed random variable. Then, the system
\eqref{sch-system8} is exactly controllable in
average with $H=L^2(G)$, $V=H^2(G)\cap H_0^1(G)$
and $U=L^2(G_0)$.

  \item Let $\a(\cd)$ be a random variable with
Laplace distribution. Then, the system
\eqref{sch-system8} is exactly controllable in
average with $H=L^2(G)$, $V=H^4(G)\cap H_0^1(G)$
and $U=L^2(G_0)$.

  \item Let $\a(\cd)$ be a random variable with
standard Chi-squared distribution. Then, the
system \eqref{sch-system8} is exactly
controllable in average with $H=L^2(G)$,
$V=H^k(G)\cap H_0^1(G)$ and $U=L^2(G_0)$.
\end{itemize}
\end{proposition}

Further, let us consider the averaged null
controllability problem for the system
\eqref{sch-system8}. Due to Remark \ref{rmk2},
we only need to prove the following
observability estimate:
\begin{equation}\label{10.20-eq3}
\(\int_\O\Big|
z(0,\o;z_0)\Big|^2_{L^2(G)}d\dbP(\o)\)^{\frac{1}{2}}
\leq C\int_E \Big|\int_\O
z(x,t,\o;z_0)d\dbP(\o)\Big|_{L^2(G_0)}dt.
\end{equation}
When  $\a(\cd)$ is a normally distributed random
variable,  we have that
\begin{equation}\label{10.20-eq5}
\begin{array}{ll}\ds
\int_\O\Big|
z(0,\o;z_0)\Big|^2_{L^2(G)}d\dbP(\o)\3n&\ds =
\frac{1}{\sqrt{2\pi}}\int_{-\infty}^\infty
e^{-\frac{\a^2}{2}}\sum_{j=1}^\infty
z_{0,j}^2e^{-2i\l_j \a
T}d\a\\
\ns&\ds= \sum_{j=1}^\infty  z_{0,j}^2e^{-2\l_j^2
t^2}.
\end{array}
\end{equation}
When  $\a(\cd)$ is a random variable with  the
Cauchy distribution, we have that
\begin{equation}\label{10.20-eq6}
\begin{array}{ll}\ds
\int_\O\Big|
z(0,\o;z_0)\Big|^2_{L^2(G)}d\dbP(\o)\3n&\ds =
\int_{-\infty}^\infty
\frac{1}{\pi(1+\a^2)}\sum_{j=1}^\infty
z_{0,j}^2e^{-2i\l_j \a
T} d\a\\
\ns&\ds= \sum_{j=1}^\infty z_{0,j}^2e^{-2\l_j
T}.
\end{array}
\end{equation}

Thanks to \eqref{10.20-eq5} and
\eqref{10.20-eq5}, and similar to the proof of
Theorem \ref{th0}, one can show that the
inequality \eqref{10.20-eq3} holds when
$\a(\cd)$ is a random variable with normal
 or Cauchy distribution. Therefore,
we have the following result:

\begin{proposition}\label{con-prop3}
Under the assumptions of Theorem \ref{con-th2}
the system \eqref{sch-system8} is null
controllable in average if $\a$ is a random
variable with normal or Cauchy distribution.
\end{proposition}


\subsubsection{Control of the variance of the state of
the system \eqref{ab-system0}}


If a system is null(\resp exactly) controllable
in average, then we can drive the average of the
state to rest (\resp a given destination). It is
then natural to analyse whether one can
 control the higher order moments of the
state and, in particular, the covariance ${\rm
Cor}(y)$ of the state $y$, which is defined as
follows:
$$
{\rm Cor}(y) = \int_{\O}\(y(T,\o)-\int_{\O}
y(T,\o) d\dbP(\o)\)\otimes \(y(T,\o)-\int_{\O}
y(T,\o) d\dbP(\o)\) d\dbP(\o),
$$
where  $\otimes$ denotes the tensor product of
two elements in $H$.

${\rm Cor}(y)$ is an element in $H\otimes H$,
which is also a Hilbert space.  With the goal of
reinforcing  the averaged controllability
property studied in this paper, a first thought
could be to look for the controllability of the
covariance. But this is not suitable since the
variance measures how far the random variable is
spread out. In other words, loosely speaking,
the smaller the $|{\rm Cor}(y(T))|_{H\otimes H}$
is, the closer the state $y(T,\cd)$ is
concentrated to an element in $H$ and,
therefore, closer to the simultaneous
controllability property. A more natural goal is
to look for the null controllability of the
covariance. However, from the definition of
${\rm Cor}(\cd)$, it is trivial to see that
${\rm Cor}(y(T))=0$ if and only if
$y(T,\cd)=\int_{\O} y(T,\o) d\dbP(\o)$ with
probability one, i.e., $y(T,\o)$ is independent
of $\o$ with probability one. Nevertheless, this
is impossible if we only use parameter
independent control (see Remark
\ref{5.28-rem1}). Therefore, since we want to
make $|{\rm Cor}(y(T))|_{V\otimes V}$ as small
as possible, it is natural to study the
following optimal control problem:

\medskip

{\bf Problem (OP2)}: Minimize $|{\rm
Cor}(y(T))|_{V\otimes V}$ for
$$
\begin{array}{ll}\ds
u(\cd)\in \cU_{null}[0,T]\3n&\ds\=\Big\{u\in
L^2(0,T;U): \mbox{ the average of solutions to
\eqref{ab-system0} corresponds }\\
\ns&\ds \q \mbox{ to $u$ satisfies that
}\int_{\O} y(T,\o;y_0)d\o=0\Big\}.
\end{array}
$$
Problem (OP2) is an optimal control problem with
a terminal constraint. Following \cite{LY,YZ},
one can derive a Pontryagin type Maximum
Principle for it. However, this is beyond the
scope of this paper.


\subsection{Open problems}


There are many interesting and important (at
least for us) problems in this topic. We present
some of them here briefly:

\begin{itemize}

\item{\bf Averaged controllability problems for
general random heat and Schr\"odinger
equations.}

We have only studied the averaged
controllability problems for some very special
classes of random heat and Schr\"odinger
equations, for which the averaged dynamics could
be computed explicitly.  It would be interesting
to investigate some more general classes of
parameter dependent systems. For example, is it
\eqref{heat-system0} null controllable in
average if we take $\a$ to be a random variable
with Chi-squared distribution?

Furthermore,  the method we use to study
\eqref{heat-system0} and \eqref{sch-system0}
depends on the fact that the eigenfunctions of
$-\a A_\D$ are independent of $\o$. Thus, more
general random heat and Schr\"odinger equations
(as, for instance, \eqref{heat-system2}) where
the principal part of the PDE depends on $\o$
cannot be treated in this way.

 The method of proof employed to show the  approximate averaged controllability of \eqref{heat-system2}
is based on the use of the space-time
analyticity  properties of solutions
\eqref{heat-system3}  to derive unique
continuation properties, but it does not provide
any quantitative information.

\medskip

\item {\bf The relationship between averaged controllability
properties and the random variable}.

We have shown that different random variable
$\a(\cd)$ may lead to different controllability
property of the system \eqref{sch-system0}. It
is interesting to give a description of the
relationship between the random variable
$\a(\cd)$ and the controllability property of
the system \eqref{sch-system0}. For instance,
for what kind of random variables, the system
\eqref{sch-system0} is null  controllable in
average? Is there a random variable such that
the system \eqref{sch-system0} is neither
exactly controllable in average nor null
controllable in average?

\medskip

\item {\bf The null averaged controllability problem for
random fractional heat equations}.

We have proven that the random fractional
Schr\"odinger equations are null  controllable
in average for $\g\in
(\frac{1}{2},\frac{1}{2})$. It is more
interesting to study the same problem for random
fractional heat equations, which describe the
anomalous diffusion process. In particularly,
consider the following system:
\begin{equation}\label{fheat-system0}
\frac{}{}\left\{
\begin{array}{ll}\ds
y_t + \a A^{\g}_\D y = B u &\mbox{ in } (0,T],\\
\ns\ds y(0)=y_0.
\end{array}
\right.
\end{equation}
Here $\g>0$, $y_0\in L^2(G)$, $u\in
L^2(0,T;L^2(G_0))$ and $Bu=\chi_{G_0}u$.

It is clear that when $\g\leq \frac{1}{2}$ and
$\a(\cd)$ is a uniformly or exponentially
distributed random variables, the system
\eqref{fheat-system0} is not null  controllable
in average. However, is it possible to find some
random variable $\a(\cd)$ such that the system
\eqref{fheat-system0} is  null  controllable in
average for some $\g\in (0,\frac{1}{2}]$?

\medskip

\item {\bf Averaged controllability problems for general
random evolution partial differential
equations}.

 The method used in this paper can
only be applied to  the case of equations in
which the random operator is $\a\D$, which is
very restrictive.

Several methods have been developed  to solve
the controllability problems for deterministic
partial differential equations. For the heat
equation, the existing methods include
 Carleman estimates (\cite{FI}) and the
moment method (\cite{AI}). For Schr\"odinger
equations, Carleman estimates can also be
applied (\cite{LTZ}), together with the moment
method (\cite{AI}), multipliers  (\cite{Ma}),
microlocal analysis (\cite{Le1}), etc.

It would be interesting to generalize these
method to deal with the averaged controllability
problems for more general random heat and
Schr\"odinger equations. The paper \cite{ABGT}
contains an interesting survey on the
controllability of parabolic systems. It would
be interesting to explore possible applications
to averaged control.

\medskip

\item {\bf Averaged controllability problems for nonlinear
random evolution equations}.

We have studied the averaged controllability
problems for linear random evolution equations.
The same problem can be considered for nonlinear
random evolution equations. A possible method to
handle the nonlinear problem is to follow what
people do for the classical controllability
problems, that is, combining the controllability
result and some fixed point theorem or inverse
mapping theorem. However, as the average of the
solution of linear transport equation with
respect to velocity helps people study the
nonlinear transport equations, we expect that
one can get better result than the ones obtained
by applying the method mentioned above directly.
For example, can one prove that a random
Schr\"odinger equation with a cubic term is
exactly or null controllable in average?

\medskip

\item {\bf Numerical approximation of  averaged
controls.}

The numerical approximation of control problems
is studied extensively in the literature. We
refer the readers to \cite{EZ,GLH,Z3} and the
rich references therein for this topic. This
question also arises in the context of averaged
control. This can be done based on the
variational characterization of the average
controls given in this paper.

\end{itemize}

\appendix
\section{Appendix}


\subsection{Reduction of  averaged
controllability to averaged
observability}\label{sec-ab-con}


We have the following results.
\begin{theorem}\label{th-ex-ob}
System \eqref{ab-system0} is exactly
controllable in average in $E$ with the control
cost $C>0$ if and only if the adjoint system
\eqref{ab-system1} is exactly observable in
average in $E$.
\end{theorem}

{\it Proof of Theorem \ref{th-ex-ob}}\,: The
``if" part.  Let us define a linear subspace
$\cX\subset L^2(E;U)$ as
$$
\cX\=\Big\{\chi_E(\cd) \int_\O B^*(\o)
z(\cd,\o;z_0) d\dbP(\o):\,  z_0\in V'\Big\}.
$$
and a linear functional $F$ on $\cX$ as
$$
F\(\chi_E(\cd) \int_\O B^*(\o) z(\cd,\o;z_0)
d\dbP(\o)\)=\langle y_1, z_0\rangle_{V,V'} -
\Big\langle y_0, \int_\O
z(0,\o;z_0)d\dbP(\o)\Big\rangle_{V,V'}.
$$
From \eqref{def-ex-ob-eq1}, we know that
$$
\begin{array}{ll}\ds
\q F\(\chi_E(\cd)\int_\O B^*(\o)
z(\cd,\o;z_0)\dbP(\o)\)\\
\ns\ds\leq C\big(|y_0|_V +
|y_1|_V\big)\(|z_0|_{V'} + \Big|\int_\O
z(0,\o;z_0)d\dbP(\o)\Big|_{V'}\)\\
\ns \ds \leq C\big(|y_0|_V +
|y_1|_V\big)\Big|\chi_E(\cd)\int_\O B^*(\o)
z(\cd,\o;z_0)\dbP(\o)\Big|_{L^2(0,T;U)}.
\end{array}
$$
Hence, $F$ is a bounded linear functional on
$\cX$ with norm $|F|_{\cL(\cX,\dbR)}\leq
C(|y_0|_V + |y_1|_V)$. Then, it can be extended
to a bounded linear functional on $L^2(E;U)$
with the same norm. We still denote by $F$ the
extension if there is no confusion. Then, by
Riesz Representation Theorem, there is a
$u(\cd)\in L^2(E;U)$ such that for any
$v(\cd)\in L^2(E;U)$,
$$
F(v)=\langle v,u\rangle_{L^2(E;U)}
$$
and
$$
|u|_{L^2(E;U)}=|F|_{\cL(\cX,\dbR)}\leq C(|y_0|_V
+ |y_1|_V).
$$
From the definition of $F(\cd)$, we know that
for any $z_0\in V'$,
\begin{equation}\label{6.9-eq3}
\langle y_1, z_0\rangle_{V,V'} \!- \Big\langle
y_0, \int_\O\!
z(0,\o;z_0)d\dbP(\o)\Big\rangle_{V,V'}\! \!=\!
\int_0^T\!\!\chi_E(t) \Big\langle u(t), \int_\O
\!B^*(\o) z(t,\o;z_0) d\dbP(\o)\Big\rangle_U \!
dt.
\end{equation}
We claim that $u(\cd)$ is the control we need.
Indeed, taking the dual product of $V'$ and $V$
of $z=z(t,\o;z_0)$ with \eqref{ab-system0} and
integrating the product with respect to $t$ in
$(0, T)$ and $\o$ in $\O$,  we obtain that
\begin{equation}\label{6.9-eq4}
\begin{array}{lll}
\ds\q\int_0^T \chi_E(t)\Big\langle u(t), \int_\O B^*(\o)z(t,\o;z_0) d \dbP(\o)\Big\rangle_U  dt \\
\ns\ds
=\int_0^T \chi_E(t)\int_\O \langle u(t), B^*(\o) z(t,\o;z_0)\rangle_U d\dbP(\o) dt \\
\ns\ds=\int_0^T\chi_E(t) \int_\O \langle B(\o) u(t), z(t,\o;z_0)\rangle_{V,V'} d\dbP(\o) dt \\
\ns\ds= \int_\O \langle y(T,\o;y_0), z_0\rangle_{V,V'}d\dbP(\o) - \int_\O \langle y_0, z(0,\o;z_0)\rangle_{V,V'}d\dbP(\o)\\
\ns\ds = \Big\langle\int_\O y(T,\o;y_0)
d\dbP(\o), z_0\Big\rangle_{V,V'} - \Big\langle
y_0, \int_\O
z(0,\o;z_0)d\dbP(\o)\Big\rangle_{V,V'}.
\end{array}
\end{equation}
From \eqref{6.9-eq3} and \eqref{6.9-eq4}, we
conclude that
$$
\Big\langle\int_\O y(T,\o;y_0) d\dbP(\o),
z_0\Big\rangle_{V,V'}=\langle y_1,
z_0\rangle_{V,V'},\q \forall\,z_0\in V',
$$
which deduces that $ \int_\O y(T,\o;y_0)
d\dbP(\o)=y_1$.

\vspace{0.2cm}

The ``only if" part. Let $z_0\in V'$. We choose
$y_0, y_1\in V'$ which satisfy that
$$
\left\{
\begin{array}{ll}\ds
|y_0|_{V'} = |y_1|_{V'}\leq 2,\\
\ns\ds -\Big\langle y_0, \int_\O
z(0,\o;z_0)d\dbP(\o)\Big\rangle_{V,V'} =
\Big|\int_\O
z(0,\o;z_0)d\dbP(\o)\Big|_{V'}^2,\\
\ns\ds \langle y_1, z_0\rangle_{V,V'} =
|z_0|_{V'}^2.
\end{array}
\right.
$$
Let $u(\cd)$ be the control such that
\begin{equation}\label{6.9-eq5}
|u|_{L^2(E;U)}\leq C(|y_0|_V + |y_1|_V)\leq C
\end{equation}
and
$$
\int_\O y(T,\o;y_0)d\dbP(\o)=y_1.
$$
Then, from \eqref{6.9-eq4}, we have that
$$
|z_0|_{V'} + \Big|\int_\O
z(0,\o;z_0)d\dbP(\o)\Big|_{V'} =
\int_0^T\chi_E(t) \Big\langle u(t), \int_\O
B^*(\o) z(t,\o;z_0) d\dbP(\o)\Big\rangle_U dt.
$$
Thus, we find that
$$
|z_0|_{V'} + \Big|\int_\O
z(0,\o;z_0)d\dbP(\o)\Big|_{V'}\leq
|u|_{L^2(E;U)}\Big|\int_\O B^*(\o) z(\cd,\o;z_0)
d\dbP(\o)\Big|_{L^2(E;U)}.
$$
This, together with \eqref{6.9-eq5}, implies
that
$$
|z_0|_{V'}\leq C\Big|\int_\O B^*(\o)
z(\cd,\o;z_0) d\dbP(\o)\Big|_{L^2(E;U)}.
$$
\begin{theorem}\label{th-nu-ob}
System \eqref{ab-system0} is null controllable
in average in $E$ with control the cost $C>0$ if
and only the adjoint system \eqref{ab-system1}
is null observable in average.
\end{theorem}
The proof of Theorem \ref{th-nu-ob} is very
similar to the one for Theorem \ref{th-ex-ob}.
We omit it here.

\begin{theorem}\label{th-ap-ob}
System \eqref{ab-system0} is approximately
controllable in average in $E$  if and only if
the adjoint system \eqref{ab-system1} satisfies
the averaged unique continuation property in
$E$.
\end{theorem}
{\it Proof of Theorem \ref{th-ap-ob}}\,: The
``if" part. Since the system \eqref{ab-system0}
is linear, we may assume that $y_0=0$. Then, we
only need to prove the following set
$$
\cA_T\=\Big\{\int_\O y(T,\o;0)d\dbP(\o):\, y
\mbox{ solves \eqref{ab-system0} with some
control }u(\cd) \Big\}
$$
is dense in $V$.  We do this by contradiction
argument. If $\cA_T$ was not dense in $V$, then
we can find a $\f\in V'$ with $|\f|_{V'}=1$ such
that
$$
\langle \psi,\f
\rangle_{V,V'}=0,\;\forall\,\psi\in\cA_T.
$$
On the other hand, similar to \eqref{6.9-eq4},
we have that
\begin{equation}\label{6.9-eq6}
\int_0^T \chi_E(t)\Big\langle u(t), \int_\O
B^*(\o)z(t,\o;z_0) d\dbP(\o)\Big\rangle_U dt =
\Big\langle\int_\O y(T,\o;0) d\dbP(\o),
z_0\Big\rangle_{V,V'}.
\end{equation}
Let $z_0=\f$ in \eqref{6.9-eq6}. We have that
$$
\int_0^T \chi_E(t)\Big\langle u(t), \int_\O
B^*(\o)z(t,\o;\f) d\dbP(\o)\Big\rangle_U
dt=0,\q\forall\, u(\cd)\in L^2(0,T;U).
$$
Hence, we find that
$$
\chi_E(\cd)\int_\O B^*(\o)z(\cd,\o;\f)
d\dbP(\o)=0 \q\mbox{ in } L^2(0,T;U),
$$
which implies that $\f=0$ and leads to a
contradiction.

The ``only if" part. We utilize the
contradiction argument again. We assume that
there is a $z_0\in V'$ with $|z_0|_{V'}=1$, such
that
$$
\chi_E(\cd)\int_\O B^*(\o)z(\cd,\o;z_0)
d\dbP(\o)=0 \q\mbox{ in } L^2(0,T;U).
$$
This, together with \eqref{6.9-eq6}, implies
that
\begin{equation}\label{6.9-eq7}
\Big\langle\int_\O y(T,\o;0) d\dbP(\o),
z_0\Big\rangle_{V,V'}=0,\q \forall\, u(\cd)\in
L^2(E;U).
\end{equation}
On the other hand, since \eqref{ab-system0} is
approximately controllable in average, we can
find a $u(\cd)\in L^2(E;U)$ such that
$$
\Big|\int_\O y(T,\o;0) d\dbP(\o)-
z_0\Big|_{V,V'}<\frac{1}{2}.
$$
It is clear that for this $\int_\O y(T,\o;0)
d\dbP(\o)$,
$$
\Big\langle\int_\O y(T,\o;0) d\dbP(\o),
z_0\Big\rangle_{V,V'}>\frac{1}{2},
$$
which contradicts \eqref{6.9-eq7}.
\endpf

\begin{remark}\label{rmk2}
One can also consider the case where $y_0\in
L^2(\O;V)$, i.e. when the datum to be controlled
depends on the parameter $\o$ as well.

The proof of Theorem \ref{th-ex-ob} applies
replacing the terms $ \langle y_0, \int_\O
z(0,\o;z_0)d\dbP(\o) \rangle_{V,V'}$ and
$|y_0|_V$ by the terms $\int_\O \langle y_0(\o),
z(0,\o;z_0) \rangle_{V,V'}d\dbP(\o)$ and
$|y_0(\cd)|_{L^2(\O;V)}$, respectively. The same
can be said about  Theorem \ref{th-ap-ob}.

The situation is different and much more
delicate in the context of averaged null
controllability. Note that, when considering
initial data to be controlled depending on $\o$,
one requires an observability estimate of the
form
\begin{equation}\label{10.20-eq4}
\int_\O \Big|z(0,\o;z_0)\Big|_{V'}^2d\dbP(\o)
\leq C\int_0^T \chi_E(t)\Big|\int_\O
B(\o)^*z(t,\o;z_0)d\dbP(\o)\Big|_{U}^2dt,
\end{equation}
which, in principle, is much stronger than the
one we have in which, we get observability
estimates on
$$
\Big|\int_\O z(0,\o;z_0)d\dbP(\o)\Big|_{V'}^2.
$$

This difficulty does not arise in the context of
exact averaged controllability where, we recover
the norm of $z_0$, and this yields also
estimates on $z(0,\o;z_0)$ for all $\o$ and in
particular on  $\int_\O
\Big|z(0,\o;z_0)\Big|_{V'}^2d\dbP(\o)$.

But the property of null averaged controllability with initial data independent of $\o$ does not seem to suffice to derive the same property with initial data that depend on $\o$. This is an interesting issue for further work.%
\end{remark}
%


\subsection{Variational characterization of
the controls of minimal norm}\label{sec-var}


\q\;We have shown the existence of the exactly
averaged control(\resp null averaged control,
approximately averaged control), provided that
the adjoint system is exactly averaged
observable(\resp null averaged observable,
satisfying averaged unique continuation
property). These results allow concluding
whether a system is controllable in an averaged
sense.  In this section, we give variational
characterizations of the controls. Such kind of
results not only derive characterizations of the
controls but also serve for computational
purposes.

\begin{theorem}\label{th-ex-va}
If the system \eqref{ab-system1} is exactly
averaged observable, then the exact averaged
control for the system \eqref{ab-system0} of
minimal $L^2(0,T; U)$-norm is given by
\begin{equation}\label{th-ex-va-eq1}
u(t) = \chi_E(t)\int_\O B^*(\o) z(0,\o; \hat
z_0) \dbP(\o),
\end{equation}
where   $\hat z_0\in V'$ minimizes the
functional
\begin{equation}\label{th-ex-va-eq2}
\begin{array}{lll}\ds
J\left(z_0\right) \3n&\ds= \frac{1}{2} \int^T_0
\chi_E(t)\Big| \int_\O\!\! B^*(\o) z(t,\o;z_0) d
\dbP(\o)\Big|^2_{U} dt  -
\langle y_1, z_0\rangle_{V,V'}\\
\ns&\ds\q + \Big\langle y_0, \int_\O
z(0,\o;z_0)d\dbP(\o)\Big\rangle_{V,V'}.
\end{array}
\end{equation}
\end{theorem}
\begin{remark}
The control given by \eqref{th-ex-va-eq1} is an
average of functions of the form $B(\cd)^* \hat
z(t,\cd)$. For each sample point $\o$, $B(\o)^*
\hat z(t,\o)$ can be chosen to be a control.
However, generally speaking, it does not steer
the initial datum $y_0$ to the final one $y_1$.
\end{remark}
{\it Proof of Theorem \ref{th-ex-va}}\,: Define
a functional $J(\cd):V'\to\dbR$ as follows:
\begin{equation}\label{6.10-eq2}
\begin{array}{ll}\ds
J(z_0) \3n&\ds= \frac{1}{2}\int^T_0
\chi_E(t)\Big| \int_\O B^*(\o)z(t,\o;z_0) d
\dbP(\o)\Big|^2_U dt - \langle y_1,
z_0\rangle_{V,V'}\\
\ns&\ds \q + \Big\langle y_0, \int_\O
z(0,\o;z_0)d\dbP(\o)\Big\rangle_{V,V'}.
\end{array}
\end{equation}

It is easy to see that the functional $J:
V'\to\dbR$ is continuous and convex. From
\eqref{def-ex-ob-eq1}, we have that $J(\cd)$ is
coercive. Then, we know that $J(\cd)$ has a
unique minimizer $\hat z_0$. Let $\hat
z(\cd,\cd)$ be the corresponding solution of the
adjoint system. By computing the first variation
of $J(\cd)$, it can be seen that
\begin{equation}\label{6.10-eq1}
\begin{array}{ll}\ds
\langle y_1, z_0\rangle_{V,V'} \3n&\ds =
\int_0^T \chi_E(t)\Big\langle \int_\O B^*(\o)
z(t, \o; \hat z_0) d \o, \int_\O B^*(\o)
z(t,\o;z_0)
d\dbP(\o)\Big\rangle_U  dt\\
\ns&\ds \q  +\Big\langle y_0, \int_\O
z(0,\o;z_0)d\dbP(\o)\Big\rangle_{V,V'}, \,
\forall z_0 \in V'.
\end{array}
\end{equation}
From \eqref{6.10-eq1}, we know that if we choose
the control as \eqref{th-ex-va-eq1}, then
$$
\int_\O y(T,\o;y_0)d\dbP(\o)=y_1.
$$

Now we prove that $u(\cd)$ given by
\eqref{th-ex-va-eq1} is the control with minimal
$L^2(0,T;U)$-norm, which drives the mathematical
expectation of the solution of the system
\eqref{ab-system0} from $y_0$ to $y_1$. Let us
choose $z_0=\hat z_0$ in \eqref{6.10-eq1}. We
have that
\begin{equation}\label{6.10-eq3}
\langle y_1, \hat z_0\rangle_{V,V'} -
\Big\langle y_0, \int_\O z(0,\o;\hat
z_0)d\dbP(\o)\Big\rangle_{V,V'} =
\int_0^T\chi_E(t) \Big|\int_\O B^*(\o) z(t, \o;
\hat z_0) d \dbP(\o) \Big|^2_U  dt.
\end{equation}
Let $\tilde u(\cd)$ be another control which
steers the mathematical expectation of the
solution to the system \eqref{ab-system0} from
$y_0$ to $y_1$. Then we obtain that
\begin{equation}\label{6.10-eq4}
\begin{array}{ll}\ds
\langle y_1, \hat z_0\rangle_{V,V'} \3n&\ds =
\int_0^T \chi_E(t)\Big\langle \tilde u(t),
\int_\O B^*(\o) z(t,\o;z_0)
d\dbP(\o)\Big\rangle_U  dt\\
\ns&\ds \q  +\Big\langle y_0, \int_\O z(0,\o;
\hat z_0)d\dbP(\o)\Big\rangle_{V,V'},  \,
\forall z_0 \in H.
\end{array}
\end{equation}
From \eqref{6.10-eq3} and \eqref{6.10-eq4}, we
get that
$$
\begin{array}{ll}\ds
\q\int_0^T \chi_E(t)\Big|\int_\O B^*(\o) z(0,\o;
\hat z_0) d \dbP(\o) \Big|^2_U  dt
\\
\ns\ds =\int_0^T \chi_E(t)\Big\langle \tilde
u(t), \int_\O B^*(\o) z(0,\o; \hat z_0)
d\dbP(\o)\Big\rangle_U dt
\\
\ns \ds \leq \Big|\chi_E \int_\O B^*(\o) z(0,\o;
\hat z_0) d\dbP(\o) \Big|_{L^2(0,T;U)}|\chi_E
\tilde u|_{L^2(0,T;U)},
\end{array}
$$
which implies that
$$
\Big|\chi_E\int_\O B^*(\o) z(0,\o; \hat z_0)
d\dbP(\o) \Big|_{L^2(0,T;U)}\leq |\chi_E\tilde
u|_{L^2(0,T;U)}.
$$
\endpf

\begin{remark}
The control is the solution of the adjoint
system with the final datum which is the
minimizer of the quadratic, convex and coercive
functional $J(\cd)$ in $V'$. One can use
numerical  methods such as implementing gradient
like iterative algorithms to compute it (see
\cite{EZ} for example). However, one will meet
the same difficulty as employing this method to
compute the control for the exact control
problems of PDEs.
\end{remark}

Similar to Theorem \ref{th-ex-va}, we can prove
the following result.
\begin{theorem}\label{th-nu-va}
If the system \eqref{ab-system1} is null
averaged observable, then the null averaged
control of minimal $L^2(0,T; U)$-norm is given
by
\begin{equation}\label{th-nu-va-eq1}
u(t) = \chi_E(t)\int_\O B^*(\o)z(0,\o; \hat z_0)
d\dbP(\o),
\end{equation}
where   $z_0\in V'$  minimizes the functional
\begin{equation}\label{th-nu-va-eq2}
\begin{array}{lll}\ds
J\left(z_0\right)=\frac{1}{2}
\int^T_0\chi_E(t)\Big| \int_\O B^*(\o)
z(t,\o;z_0) d \dbP(\o)\Big|^2_U dt  +
\Big\langle y_0, \int_\O
z(0,\o;z_0)d\dbP(\o)\Big\rangle_{V,V'}.
\end{array}\end{equation}
\end{theorem}
\begin{theorem}\label{th-ap-va}
Suppose that the system \eqref{ab-system1}
satisfies the averaged unique continuation
property. Then for any $\e>0$, the approximately
averaged control  is given by
\begin{equation}\label{th-ap-va-eq1}
u(t) = \chi_E(t)\int_\O B^*(\o) \hat z(t,
\o)d\dbP(\o).
\end{equation}
Here $\hat z$ is the solution to the adjoint
system \eqref{ab-system1} corresponding to the
datum $z_0\in V'$ which minimizes the functional
\begin{equation}\label{th-ap-va-eq3}
\begin{array}{ll}\ds
J_\e(z_0)\3n&\ds=\frac{1}{2} \int^T_0\chi_E(t)
\Big| \int_\O B^*(\o) z(t,\o;z_0) d
\dbP(\o)\Big|^2_U dt - \langle y_1,
z_0\rangle_{V,V'}\\
\ns&\ds \q + \Big\langle y_0, \int_\O
z(0,\o;z_0)d\dbP(\o)\Big\rangle_{V,V'}  + \e
|z_0|_{V'}.
\end{array}
\end{equation}
\end{theorem}

Borrowing some idea in \cite{Z1}, we can prove
the following result stronger than Theorem
\ref{th-ap-va}.

\begin{theorem}\label{th-ap-va1}
Suppose that the system \eqref{ab-system1}
satisfies the averaged unique continuation
property. For any $\e>0$ and any finite
dimensional space $X\subset V$, the solution to
\eqref{ab-system0} with the control
\begin{equation}\label{th-ap-va1-eq1}
u(t) = \chi_E(t)\int_\O B^*(\o) \hat z(t,
\o)d\dbP(\o)
\end{equation}
satisfies that
\begin{equation}\label{th-ap-va1-eq2}
\Big|y_1-\int_\O
y(T,\o;y_0)d\dbP(\o)\Big|_V<\e,\q \Pi_X y_1 =
\Pi_X \int_\O y(T,\o;y_0)d\dbP(\o).
\end{equation}
Here $\Pi_X$ denotes the orthogonal projection
from $V$  to $X$ and $\hat z(\cd)$ is the
solution to the adjoint system
\eqref{ab-system1} corresponding to the final
datum $\hat z_0\in V'$ which minimizes the
functional
\begin{equation}\label{th-ap-va1-eq3}
\begin{array}{lll}\ds
J_\e(z_0)\3n&\ds= \frac{1}{2}
\int^T_0\chi_E(t)\Big| \int_\O B^*(\o)
z(t,\o;z_0) d \dbP(\o)\Big|^2_U dt  - \langle
y_1,
z_0\rangle_{V,V'}\\
\ns&\ds \q + \Big\langle y_0, \int_\O
z(0,\o;z_0)d\dbP(\o)\Big\rangle_{V,V'} + \e
|(I-\Pi_X^*)z_0|_{V'}.
\end{array}\end{equation}
\end{theorem}

{\it Proof}\,: Clearly, $J_\e(\cd)$ is
continuous and convex. We only need to show that
it is coercive. For this, we prove that
\begin{equation}\label{6.10-eq5}
J_\e(z_0) \to \infty \mbox{ as } |z_0|_{V'}\to
\infty.
\end{equation}
Let $\{z_0^j\}_{j=1}^\infty\subset V'$ be a
sequence such that $|z_0^j|_{V'}\to\infty$ as
$j\to\infty$. Put $\check
z_0^j=|z_0^j|_{V'}^{-1}z_0^j$ for $j\in \dbN$.
Then
\begin{equation}\label{6.10-eq6}
\begin{array}{ll}\ds
\frac{z_0^j}{|z_0^j|_{V'}}\3n&\ds=\frac{1}{2}|z_0^j|_{V'}
\int_0^T\chi_E(t)\Big|\int_\O
B(\o)^*z(t,\o;\check z_0^j)d\dbP(\o)\Big|_U^2 dt
- \langle y_1,
\check z_0^j\rangle_{V,V'}\\
\ns&\ds \q + \Big\langle y_0, \int_\O
z(0,\o;\check z_0^j)d\dbP(\o)\Big\rangle_{V,V'}
+ \e|(I-\Pi_X^*)\check z_0^j|_{V'}.
\end{array}
\end{equation}
If
$$
\liminf_{j\to\infty}\int_0^T\chi_E(t)\Big|\int_\O
B(\o)^* z(t,\o;\check
z_0^j)d\dbP(\o)\Big|_U^2dt>0,
$$
then we get from \eqref{6.10-eq6} that
$$
J_\e(z_0^j)\to\infty \;\mbox{ as }\; j\to\infty.
$$
Hence, we only need to consider the case that
$$
\liminf_{j\to\infty}\int_0^T\chi_E(t)\Big|\int_\O
B(\o)^* z(t,\o;\check
z_0^j)d\dbP(\o)\Big|_U^2dt=0.
$$
Since $\{\check z_0^j\}_{j=1}^\infty$ is bounded
in $V'$, we can find a subsequence of it, which
is still denoted by $\{\check
z_0^j\}_{j=1}^\infty$ if there is no confusion,
such that $\check z_0^j$ converges to some
$\check z_0$ in $V'$ weakly. Thus, we have that
$ \chi_E(\cd)\!\int_\O\! B(\o)^*z(\cd,\o;\check
z_0^j)d\dbP(\o)$
 converges to $\chi_E(\cd)\!\int_\O\!
B(\o)^*z(\cd,\o;\check z_0)d\dbP(\o)$ weakly in
$L^2(0,T;U)$. Then, we have that
$$
\int_0^T\chi_E(t)\Big|\int_{\O}B(\o)^*z(t,\o;\check
z_0)d\dbP(\o)\Big|_U^2 dt \leq
\liminf_{j\to\infty}\int_0^T\chi_E(\cd)\Big|
\int_{\O}B(\o)^*z(t,\o;\check
z_0)d\dbP(\o)\Big|_U^2 dt=0,
$$
which implies that $\chi_E
\int_{\O}B(\o)^*z(t,\o;\check z_0)d\dbP(\o)=0$
in $L^2(0,T;U)$. Thus, we know that $\check
z_0=0$. Since $X$ is finite dimensional, we have
that
$$
\lim_{j\to\infty} |(I-\Pi_X^*)\check z_0^j|_{V'}
= 1.
$$
Therefore,
$$
\liminf_{j\to\infty}
\frac{J_\e(z_0^j)}{|z_0^j|_{V'}} \geq
\liminf_{j\to\infty} \( - \langle y_1, \check
z_0^j\rangle_{V,V'}\! +\! \Big\langle y_0,
\int_\O z(0,\o;\check
z_0^j)d\dbP(\o)\Big\rangle_{V,V'}\! +\! \e
\)=\e,
$$
which implies that $J_\e(z_0^j)\to\infty$ as
$|z_0^j|_{V'}\to\infty$. By this, we get the
coercivity of $J_\e(\cd)$. Hence, we know that
there is a minimizer $\hat z_0$ of $J_\e(\cd)$.

For any $\d>0$ and $z_0\in V'$, we have that
$$
0\leq \frac{1}{h}\big[ J_\e(\hat z_0 + \d z_0) -
J_\e(\hat z_0)\big],
$$
which implies that
$$
\begin{array}{ll}\ds
-\e|(I-\Pi_X^*)z_0|_{V'} \3n&\ds\leq
\int_0^T\chi_E(t) \Big\langle\int_\O
B(\o)^*z(t,\o;\hat z_0)d\dbP(\o), \int_\O
B(\o)^*z(t,\o;
z_0)d\dbP(\o) \Big\rangle_U dt \\
\ns&\ds \q - \langle y_1, z_0\rangle_{V,V'} +
\Big\langle y_0, \int_\O z(0,\o;
z_0)d\dbP(\o)\Big\rangle_{V,V'}.
\end{array}
$$
Similarly, we have that
$$
\begin{array}{ll}\ds
\e|(I-\Pi_X^*)z_0|_{V'} \3n&\ds\geq \int_0^T
\chi_E(t)\Big\langle\int_\O B(\o)^*z(t,\o;\hat
z_0)d\dbP(\o), \int_\O B(\o)^*z(t,\o;
z_0)d\dbP(\o) \Big\rangle_U dt \\
\ns&\ds \q - \langle y_1, z_0\rangle_{V,V'} +
\Big\langle y_0, \int_\O z(0,\o;
z_0)d\dbP(\o)\Big\rangle_{V,V'}.
\end{array}
$$
Hence, for any $z_0\in V'$,
$$
\begin{array}{ll}\ds
\e|(I-\Pi_X^*)z_0|_{V'} \3n&\ds\geq
\Big|\int_0^T\chi_E(t) \Big\langle\int_\O
B(\o)^*z(t,\o;\hat z_0)d\dbP(\o), \int_\O
B(\o)^*z(t,\o;
z_0)d\dbP(\o) \Big\rangle_U dt \\
\ns&\ds \q - \langle y_1, z_0\rangle_{V,V'} +
\Big\langle y_0, \int_\O z(0,\o;
z_0)d\dbP(\o)\Big\rangle_{V,V'}\Big|.
\end{array}
$$
This, together with \eqref{6.9-eq4}, implies
that for any $z_0\in V'$,
$$
\Big|\Big\langle z_0, y_1-\int_\O
y(T,\o;y_0)\Big\rangle_{V,V'}\Big| \leq
\e|(I-\Pi_X^*)z_0|_{V'}.
$$
Hence, we get \eqref{th-ap-va1-eq2}.

\section*{Acknowledgement}
This work is supported by  the Advanced Grants
NUMERIWAVES/FP7-246775 of the European Research
Council Executive Agency, FA9550-14-1-0214 of
the EOARD-AFOSR, PI2010-04 and the BERC
2014-2017 program of the Basque Government, the
 MTM2011-29306-C02-00 and
SEV-2013-0323 Grants of the MINECO a Humboldt
Research Award (University of
Erlangen-N\"urnberg), and the NSF of China under
grant 11101070.

The authors thank the CIMI - Toulouse for the
hospitality and support during the preparation
of this work in the context of the Excellence
Chair in ``PDE, Control and Numerics".

The author acknowledges Martin Lazar  for
helpful comments.

\end{document}